\documentclass[11pt]{article}

\usepackage{fullpage}
\usepackage{authblk}
\usepackage{natbib}
\usepackage{algorithm}
\usepackage{algpseudocode}
\usepackage{amsmath,amsthm,amsfonts}
\usepackage{amsmath}
\usepackage[subnum]{cases}
\usepackage{tcolorbox}
\usepackage{mathrsfs}
\usepackage{amssymb}
\usepackage{color}
\usepackage{mathrsfs}
\usepackage{empheq}
\usepackage{enumitem}
\usepackage{bm}
\usepackage{multirow}
\usepackage{booktabs}
\usepackage{makecell}
\usepackage{graphicx}
\usepackage{subcaption}
\usepackage{comment}
\usepackage{cases}
\usepackage{appendix}
\usepackage{tikz}
\usetikzlibrary{arrows,shapes}
\usepackage[colorlinks= true, linkcolor=red, citecolor=blue, urlcolor=black]{hyperref}
\usepackage{cleveref}
\usepackage{marginnote}
\usepackage{diagbox} 

\numberwithin{equation}{section}

\theoremstyle{definition}

\newtheorem{theorem}{Theorem}[section]
\newtheorem{lemma}[theorem]{Lemma}

\newtheorem{defn}[theorem]{Definition}

\DeclareMathOperator*{\argmin}{arg\,min}

\usepackage{booktabs}
\usepackage{tikz}
\usepackage{pgfplots}
\usetikzlibrary{calc,intersections}

\usetikzlibrary{shapes.geometric, arrows, positioning}

\tikzset{
	mybox/.style  = {draw, rectangle, minimum width=4cm, minimum height=0.8cm, text centered, text width=4.4cm,   
		font=\normalsize},
	box/.style  = {draw, rectangle, minimum width=2.0cm, minimum height=0.6cm, text centered, text width=3.0cm,   
		font=\normalsize},
	myarrow/.style = {line width=0.2pt, draw=black, -triangle 60, postaction={draw, line width=0.2pt, shorten >=10pt,-}}
}

\tikzstyle{arrow} = [->, >=stealth, -triangle 60]

\allowdisplaybreaks

\makeatletter
\newcommand{\leqnomode}{\tagsleft@true}
\newcommand{\reqnomode}{\tagsleft@false}
\makeatother

\begin{document}

\title{Understanding the ADMM Algorithm via High-Resolution Differential Equations}
\author[1,2]{Bowen Li}
\author[1,2]{Bin Shi\thanks{Corresponding author: \url{shibin@lsec.cc.ac.cn} } }
\affil[1]{Academy of Mathematics and Systems Science, Chinese Academy of Sciences, Beijing 100190, China}
\affil[2]{School of Mathematical Sciences, University of Chinese Academy of Sciences, Beijing 100049, China}
\date\today

\maketitle

\begin{abstract}
In the fields of statistics, machine learning, image science, and related areas, there is an increasing demand for decentralized collection or storage of large-scale datasets, as well as distributed solution methods. To tackle this challenge, the~\textit{alternating direction method of multipliers}~(\texttt{ADMM}) has emerged as a widely used approach, particularly well-suited to distributed convex optimization. However, the iterative behavior of~\texttt{ADMM} has not been well understood. In this paper, we employ dimensional analysis to derive a system of high-resolution ordinary differential equations (ODEs) for~\texttt{ADMM}. This system captures an important characteristic of~\texttt{ADMM}, called the $\lambda$-correction,  which causes the trajectory of~\texttt{ADMM} to deviate from the constrained hyperplane. To explore the convergence behavior of the system of high-resolution ODEs, we utilize Lyapunov analysis and extend our findings to the discrete~\texttt{ADMM} algorithm.  Through this analysis, we identify that the numerical error resulting from the implicit scheme is a crucial factor that affects the convergence rate and monotonicity in the discrete~\texttt{ADMM} algorithm. In addition, we further discover that if one component of the objective function is assumed to be strongly convex,  the iterative average of~\texttt{ADMM} converges strongly with a rate $O(1/N)$, where $N$ is the number of iterations.


\end{abstract}

%

\section{Introduction}
\label{sec: intro}

Since the beginning of the new century, gradient-based optimization has experienced a remarkable resurgence in the rapidly evolving field of machine learning.  This resurgence can be attributed to its computational efficiency and low memory requirements, which make it highly suitable for addressing large-scale problems. As a result, gradient-based optimization algorithms have regained their dominance and are widely utilized in various applications. The study and development of these algorithms have made a comeback, attracting considerable attention and becoming a focal point of research in the field.  

To provide a clear demonstration of our motivation, it is indeed helpful for us to start with the simplest form of convex optimization problems,  known as the unconstrained problem, which can be described as
\begin{equation}
\label{eqn: unconstrain-opt}
\min_{x \in \mathbb{R}^d} f(x).
\end{equation}
When the objective function $f$ is assumed to be smooth, one of the earliest gradient-based algorithms that comes to mind is the vanilla~\textit{gradient descent}, which has been around since the Euler era. In practical terms, the smooth problem~\eqref{eqn: unconstrain-opt} corresponds to two broad classes of statistics and inverse problems: \textit{linear regression} and \textit{linear inverse problem}. These problems share the same~\textit{least-square} form, which is commonly written as
\[
\min_{x \in \mathbb{R}^d} \|Ax - b\|^2,
\]
where $A \in \mathbb{R}^{m \times d}$ is an $m \times d$ matrix and $b \in \mathbb{R}^m$ is an $m$-dimensional vector.\footnote{Throughout this paper, the notation $\|\cdot\|$ specifically refers to the $\ell_2$-norm or the  Euclidean norm, $\|\cdot\|_2$. It is worth noting that the subscript $2$ is often omitted without any specific statement.} The vanilla~\textit{gradient descent} works effectively on the~\textit{least-square} problem and performs very well in practice. 

Since the eighties of the last century, the $\ell_1$-norm has been widely recognized as a valuable tool for characterizing sparse structures in geophysics. In the context of regression models,  the $\ell_1$-regularizer is often employed to enhance the prediction accuracy and interpretability. This is achieved by augmenting the~\textit{least-square} form with the $\ell_1$-regularizer, resulting in the objective function described as
\[
\min_{x \in \mathbb{R}^d} \Phi(x):=\|Ax - b\|^2 + \lambda\|x\|_1,
\]
where the regularization parameter $\lambda > 0$ is a tradeoff between fidelity to the measurements and noise sensitivity. In statistics, this formulation is known as~\textit{least absolute shrinkage and selection operator} (\texttt{Lasso}), which was rediscovered independently and popularized by~\citet{tibshirani1996regression}. From an optimization perspective, \texttt{Lasso} extends the smooth problem~\eqref{eqn: unconstrain-opt} to a composite problem expressed as
\begin{equation}
\label{eqn: composite-opt}
\min_{x \in \mathbb{R}^d} \Phi(x):= f(x) + g(x),
\end{equation}
where $f$ satisfies the same assumption as in~\eqref{eqn: unconstrain-opt} while $g$ is only required to be a convex function without any smoothness requirement. To solve this composite problem~\eqref{eqn: composite-opt}, the~\textit{iterative shrinkage-thresholding algorithm} (\texttt{ISTA}), a proximal version of the vanilla~\textit{gradient descent},  was proposed by~\citet{daubechies2004iterative}. However, for more general applications such as~\textit{total-variation denoising} in image science~\citep{rudin1992nonlinear} and~\textit{$\ell_1$ trend filtering} in time series~\citep{kim2009ell_1}, the commonly used form is the generalized~\texttt{Lasso}, which can be written as
\[
\min_{x \in \mathbb{R}^d} \Phi(x):=\|Ax - b\|^2 + \lambda\|Fx\|_1,
\]
where $F \in \mathbb{R}^{n \times d}$ is an $n \times d$ matrix. It is worth noting that~\texttt{ISTA} is not applicable for the generalized~\texttt{Lasso} since the proximal operator does not work when the vector $x$ is multiplied by a matrix $F$. 

A well-suite optimization method for solving the generalized~\texttt{Lasso} is the~\textit{Alternating Direction Method of Multipliers} (\texttt{ADMM}), which was first introduced in the seventies of the last century by~\citet{glowinski1975approximation} and~\citet{gabay1976dual}.  In recent years, \texttt{ADMM} has gained popularity, especially in the field of machine learning, due to its effectiveness in distributed convex optimization and its ability to handle large-scale optimization problems~\citep{boyd2011distributed}. In fact, \texttt{ADMM} is not limited to solving the generalized~\textit{Lasso} but can also be applied to other optimization problems related to the $\ell_1$-norm,  such as~\textit{least absolute deviations}~\citep{boyd2004convex, hastie2009elements, wooldridge2012introductory} and~\textit{basis pursuit}~\citep{bruckstein2009sparse}.  All of the above-mentioned problems can be formulated in a general form as
\begin{equation}
\label{eqn: admm-problem}
\begin{aligned}
  \min           & \quad f(x) + g(y) \\
\mathrm{s.t.}    & \quad Fx + Gy = h
\end{aligned}
\end{equation}
where $x \in \mathbb{R}^{d_1}$ and $y \in \mathbb{R}^{d_2}$ are variables, $F \in \mathbb{R}^{m \times d_1}$ and $G \in \mathbb{R}^{m \times d_2}$ are two matrices and $h \in \mathbb{R}^{m}$ is a $m$-dimensional vector. The two objective functions, $f$ and $g$, are assumed to be convex. The~\texttt{ADMM} iterations involve forming the augmented Lagrangian
\begin{equation}
\label{eqn: augment-lagrangian}
L_{s}(x,y; \lambda) = f(x) + g(y) + \left\langle \lambda, Fx + Gy -h \right\rangle + \frac{1}{2s} \left\| Fx + Gy -h  \right\|^2 ,
\end{equation}
where $\lambda \in \mathbb{R}^m$ is a multiplier and $s>0$ is a parameter. The iteration consists of updating $x$, $y$, and the multipier $\lambda$ according to the following equations:    
\begin{subequations}
\label{eqn: admm}
\begin{empheq}[left=\empheqlbrace]{align}
         & x_{k+1} = \argmin_{x \in \mathbb{R}^d} \left\{ f(x) + \frac{1}{2s} \|Fx + Gy_k - h + s \lambda_k\|^2  \right\},                  \label{eqn: admm-x}            \\   
         & y_{k+1} = \argmin_{x \in \mathbb{R}^d} \left\{ g(y) + \frac{1}{2s} \|Fx_{k+1} +  Gy - h + s \lambda_k\|^2 \right\},              \label{eqn: admm-y}            \\
         & \lambda_{k+1} = \lambda_k + \frac1s (Fx_{k+1} + Gy_{k+1} - h).                                                                                   \label{eqn: admm-lambda}
\end{empheq}    
\end{subequations}
In general, the implicit solution in the second iteration~\eqref{eqn: admm-y} can always be obtained directly when the convex function $g$ is the $\ell_1$-norm. However, when the convex function $f$ is assumed to be $L$-smooth,\footnote{The rigorous definition of $L$-smooth functions is shown in~\Cref{sec: prelim} (\Cref{defn: L-smooth-strongly}). } the implicit solution in the first iteration~\eqref{eqn: admm-x} requires the use of a gradient-based algorithm for approximation. Therefore, for our discussion in this paper, the focus is only on scenarios where $f$ is a quadratic function for the generalized~\texttt{Lasso} or an indicator function for the~\textit{basis pursuit}. In these cases, the implicit solution can always be obtained directly. 

%
%
%
%
%
%
%
%

%
%
%
%
%
%

When taking a simple comparison of the two algorithms above, \texttt{ISTA} and \texttt{ADMM}, it is natural for us to raise the two following questions: 
\begin{tcolorbox}
\begin{itemize}
\item Why \texttt{ADMM} works on the genralized~\texttt{Lasso} and \textit{basis pursuit} problems while \texttt{ISTA} does not work? 
\item Where are the fundamental differences between the two algorithms?  
\end{itemize}
\end{tcolorbox}
\noindent In this paper, our objective is to answer these questions by employing techniques borrowed from ordinary differential equations (ODEs) and numerical analysis. Recent advancements in these techniques have proven to be successful in characterizing the convergence rates of the vanilla~\textit{gradient descent} and shedding light on the intriguing acceleration phenomenon introduced by~\citet{nesterov1983method}. Through the utilization of Lyapunov analysis and phase-space representation, a comprehensive framework based on high-resolution ODEs has been well established in the studies~\citep{shi2022understanding, shi2019acceleration, chen2022gradient, chen2022revisiting, chen2023underdamped, li2023linear}. Moreover, this framework has also been extended in~\citep{li2022linear, li2022proximal, li2023linear} to encompass the convergence rates of \texttt{ISTA} and \texttt{FISTA}s (\textit{fast iterative shrinkage-thresholding algorithms}), which are the proximal versions of \textit{Nesterov's accelerated gradient descent} methods. Several new avenues are opened up by these advancements for further research and development of optimization algorithms.




\subsection{$\lambda$-correction: a small but essential perturbation}
\label{subsec: lambda-correction}

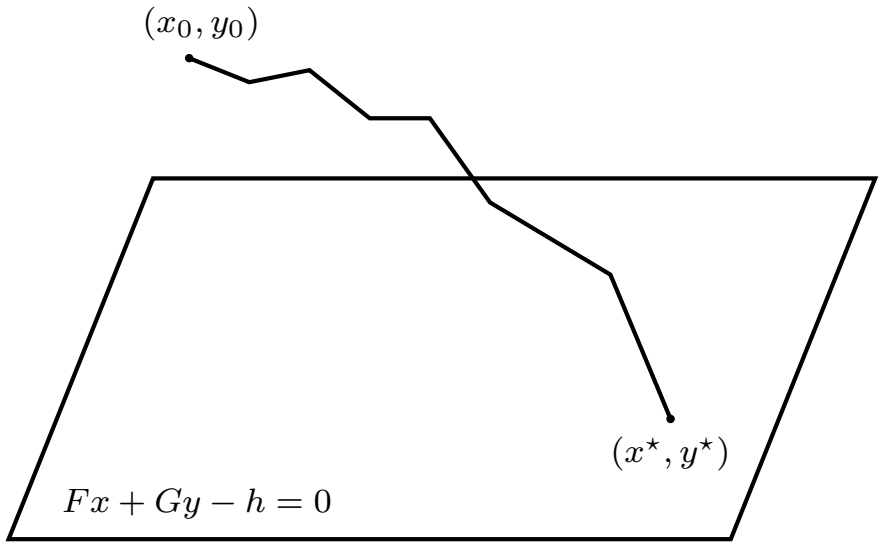
\begin{figure}[htb!]
  \centering
  \scalebox{1.6}{\begin{tikzpicture}[line width=1.0pt].
    \draw[xslant=0.4](-1,0) rectangle (5,3);
    \node at(-0.7,0.3)[right]{\scriptsize\( Fx + Gy - h = 0 \)};
%
%

    \draw[black] (0.5,4) -- (1.0, 3.8) -- (1.5,3.9)-- (2.0,3.5) -- (2.5,3.5) -- (3.0, 2.8) -- (3.5,2.5) -- (4,2.2) -- (4.5,1);
    \fill(0.5,4)circle(1pt)node[above, xshift=0.1cm]{\scriptsize\( (x_0, y_0) \)};

    \fill(4.5,1)circle(1pt)node[below]{\scriptsize\( (x^\star, y^\star) \)};

  \end{tikzpicture}}
\caption{A schematic diagram of the trajectory of~\texttt{ADMM} with any initial $(x_0, y_0)$ (Black).}
\label{fig: admm}
\end{figure}

In order to study the iterative behavior of~\texttt{ADMM}, we reformulate the third iteration~\eqref{eqn: admm-lambda} as
\begin{equation}
\label{eqn: admm-lambda-reform}
s\left(\lambda_{k+1} - \lambda_k\right) = Fx_{k+1} + Gy_{k+1} - h.
\end{equation}
From this reformulation~\eqref{eqn: admm-lambda-reform}, it can be observed that when the multiplier sequence $\{\lambda_k\}_{k=0}^{\infty}$ is introduced, except in cases where $\lambda_{k}$ is a constant starting from some integer $N$, the iterative sequence $\{(x_{k}, y_{k})\}_{k=0}^{\infty}$ generated by~\texttt{ADMM} often does not lie on the constrained hyperplane $Fx+Gy - h = 0$. In other words, the iteration of the multiplier, $s(\lambda_{k +1}- \lambda_k)$, referred to as $\lambda$-correction, acts as a small perturbation that causes the iterative sequence $\{(x_{k}, y_{k})\}_{k=0}^{\infty}$ to deviate away from the constrained hyperplane $Fx+Gy - h = 0$. This mathematical observation is visualized in~\Cref{fig: admm}.

%

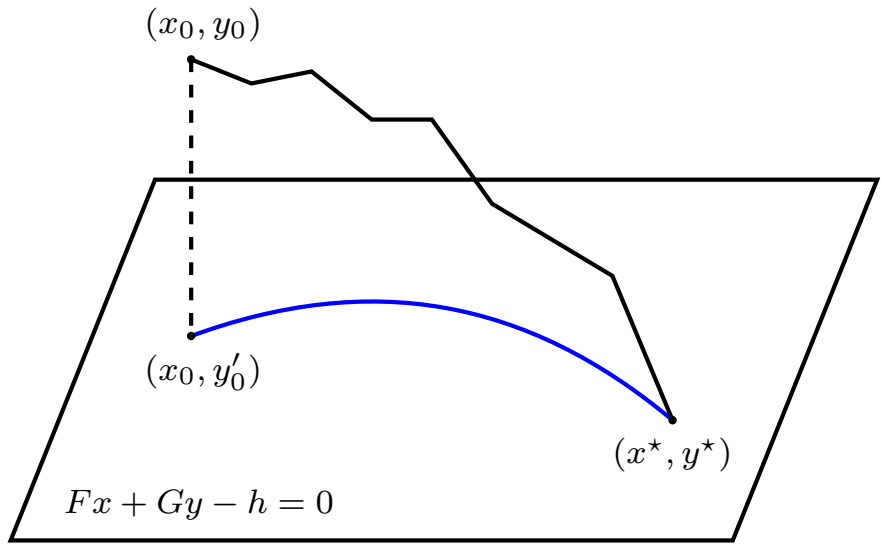
\begin{figure}[htb!]
  \centering
  \scalebox{1.6}{\begin{tikzpicture}[line width=1.0pt].
    \draw[xslant=0.4](-1,0) rectangle (5,3);
    \node at(-0.7,0.3)[right]{\scriptsize\( Fx + Gy - h = 0 \)};
  
    \path[blue] (0.5, 1.7) edge [bend left] (4.5,1);
    \fill(0.5,1.7)circle(1pt)node[below, xshift=0.1cm]{\scriptsize\( (x_0, y_0') \)};


    \draw[black] (0.5,4) -- (1.0, 3.8) -- (1.5,3.9)-- (2.0,3.5) -- (2.5,3.5) -- (3.0, 2.8) -- (3.5,2.5) -- (4,2.2) -- (4.5,1);
    \fill(0.5,4)circle(1pt)node[above, xshift=0.1cm]{\scriptsize\( (x_0, y_0) \)};

    \fill(4.5,1)circle(1pt)node[below]{\scriptsize\( (x^\star, y^\star) \)};

    \draw[black,dashed](0.5,4) -- (0.5,1.7);
  \end{tikzpicture}}
\caption{A schematic diagram of the trajectory of~\texttt{ADMM} with any initial $(x_0, y_0)$ (Black) and the continuous limit ODE with any initial $x_0$ (Blue), where the initial $y_0'$ is not arbitrary and required to satisfies $Fx_0 +Gy_0' = h$.}
\label{fig: admm-low}
\end{figure}

In recent years,  there has been a growing trend in utilizing continuous limit ODEs to model and understand discrete algorithms, particularly~\textit{Nesterov's accelerated gradient descent} methods~\citep{su2016differential, wilson2021lyapunov}. By following a similar approach, the continuous limit ODE of~\texttt{ADMM} can be expressed as  
\begin{equation}
\label{eqn: ode-low}
\dot{X} = (F^{T}F)^{-1} \left( - \nabla f(X) - F^{T} \nabla g(Y) \right). 
\end{equation}
When we set $G=I$ (identity  matrix) and $h = 0$, the continuous limit ODE~\eqref{eqn: ode-low} reduces to a generalized gradient flow as
\[
\dot{X} = (F^{T}F)^{-1} \left( - \nabla f(X) - F^{T} \nabla g(FX) \right),
\]
which was first derived in~\citep{franca2018admm}. However, an important feature is omitted in the process of taking the continuous limit.  When considering the continuous limit for the third iteration~\eqref{eqn: admm-lambda}, it is not hard for us to observe that the pair of variables $(X, Y)$ satisfy $FX + GY = h$. In other words, the pair of variables $(X, Y)$ lies on the constrained hyperplane $Fx + Gy = h$. This indicates that the continuous limit ODE~\eqref{eqn: ode-low} fails to capture the effect of $\lambda$-correction, which causes the iterative sequence $\{ (x_k, y_k) \}_{k=0}^{\infty}$ to deviate from the constrained hyperplane $Fx + Gy = h$, as depicted in~\Cref{fig: admm-low}.

%
%
%
%
%

\subsection{Numerical errors: byproduct of the implicit scheme}
\label{subsec: error-implicit}

In this context, we briefly describe an important property of a nonlinear dynamical system and its numerical discretization. Particularly, we focus on the implicit scheme, which can lead to numerical errors that have an impact on the convergence rate. We start by considering a finite-dimensional nonlinear dynamical system on $\mathbb{R}^d$, represented as
\begin{equation}
\label{eqn: non-ode-ex}
\dot{W} = \xi(W).
\end{equation}
Let $w^{\star}$ be an equilibrium, that is, $\xi(w^{\star}) = 0$. For the nonlinear system~\eqref{eqn: non-ode-ex}, we examine a Lyapunov function that has a quadratic form as  
\begin{equation}
\label{eqn: Lyapunov-ode-ex}
\mathcal{E}(t) = \frac12 \| W - w^{\star} \|^2.
\end{equation}
If the angle between the vector field $\xi(W)$ in~\eqref{eqn: non-ode-ex}  and the direction $W - w^{\star}$  is always obtuse, then the time derivative of the Lyapunov function~\eqref{eqn: Lyapunov-ode-ex} satisfies
\begin{equation}
\label{eqn: Lyapunov-ode-montone-ex}
\frac{d\mathcal{E}}{dt} = \big\langle \xi, W - w^{\star} \big\rangle \leq 0.
\end{equation}
This implies that the Lyapunov function~\eqref{eqn: Lyapunov-ode-ex} can describe that the variable $W$ converges towards the equilibrium $w^{\star}$ along the nonlinear system~\eqref{eqn: non-ode-ex}. 

The implicit scheme of the nonlinear system~\eqref{eqn: non-ode-ex} can be expressed as
\begin{equation}
\label{eqn: non-ex-implicit}
w_{k+1} - w_{k} = s\xi (w_{k+1}).
\end{equation}
To align with the continuous Lyapunov function~\eqref{eqn: Lyapunov-ode-ex}, we can write down the discrete Lyapunov function as
\begin{equation}
\label{eqn: Lyapunov-ex-implicit}
\mathcal{E}(k) = \frac12 \| w_k - w^{\star} \|^2.
\end{equation}
By utilizing the implicit scheme~\eqref{eqn: non-ex-implicit} and the discrete Lyapunov function~\eqref{eqn: Lyapunov-ex-implicit}, we can calculate the iterative difference as
\begin{equation}
\label{eqn: Lyapunov-montone-ex-implicit}
\mathcal{E}(k+1) - \mathcal{E}(k) = s \left\langle \xi(w_{k+1}), w_{k+1} - w^{\star} \right\rangle - \frac{1}{2} \|w_{k+1} - w_{k} \|^2.
\end{equation}
Similarly, if the angle between the vector field $\xi$ in~\eqref{eqn: non-ex-implicit}  and the direction $w - w^{\star}$  is always obtuse, then the iterative difference satisfies
\begin{equation}
\label{eqn: Lyapunov-montone-ex-implicit-2}
\mathcal{E}(k+1) - \mathcal{E}(k) \leq  - \frac{1}{2} \|w_{k+1} - w_{k} \|^2,
\end{equation}
where the quadratic term on the right-hand side of the inequality represents the numerical error caused by the implicit scheme~\eqref{eqn: non-ex-implicit}. Summarizing the inequalities~\eqref{eqn: Lyapunov-montone-ex-implicit-2} from $0$ to $N$, it is straightforward for us to derive the following inequality as
\begin{equation}
\label{eqn: convergence-rate-ex}
\frac{1}{2} \sum_{k=0}^{N} \|w_{k+1} - w_{k} \|^2 \leq \mathcal{E}(0).
\end{equation}
Furthermore, we can determine the convergence rates in terms of average and minimization as
\[
\frac{1}{N+1}\sum_{k=0}^{N} \|w_{k+1} - w_{k} \|^2 \leq \frac{\|w_0 - w^{\star}\|^2}{N+1} \quad \mathrm{and} \quad \min_{0 \leq k \leq N} \|w_{k+1} - w_{k} \|^2 \leq \frac{\|w_0 - w^{\star}\|^2}{N+1}. 
\]

In addition, if the angle between the vector field $\xi(w)$ in~\eqref{eqn: non-ode-ex}  and the direction $w - w^{\star}$  is not always obtuse, it implies that the time derivative of the Lyapunov function~\eqref{eqn: Lyapunov-ode-ex} is not always less than or equal to zero. In other words, it may be larger than zero at times, which means that the continuous Lyapunov function~\eqref{eqn: Lyapunov-ode-ex} does not exhibit a monotone behavior. Correspondingly, when considering the implicit scheme~\eqref{eqn: non-ex-implicit} and the discrete Lyapunov function~\eqref{eqn: Lyapunov-ex-implicit}, we can deduce that the iterative difference satisfies~\eqref{eqn: Lyapunov-montone-ex-implicit}. By maintaining the same condition from the continuous ODE to the discrete algorithm, we know that the first term on the right-hand side of the inequality~\eqref{eqn: Lyapunov-montone-ex-implicit} may occasionally be larger than zero. However, the numerical error, represented by the second term on the right-hand side of the inequality~\eqref{eqn: Lyapunov-montone-ex-implicit}, may counterbalance it in such a way that it is always nonpositive. As a result, the discrete Lyapunov function may exhibit a monotonically decreasing behavior.

\subsection{Overview of contributions}
\label{subsec: overview}

In this paper, we employ techniques borrowed from ODEs and numerical analysis to understand and analyze the~\texttt{ADMM} algorithm and highlight our contributions as follows. 

\begin{itemize}
\item[(1)] We utilize the dimensional analysis proposed in~\citep{shi2022understanding, shi2020learning, shi2021hyperparameters} to derive the system of high-resolution ODEs\footnote{The continuous limit ODE is also called the low-resolution ODE. See~\Cref{subsec: derivation-ode} for details.}  of~\texttt{ADMM} as
\begin{subequations}
\label{eqn: admm-high-ode}
\begin{empheq}[left=\empheqlbrace]{align}
        & F^{\top}G\dot{Y}  = F^{\top}\Lambda + \nabla f(X),             \label{eqn: admm-ode-x}               \\
        & 0 = G^{\top}\Lambda  + \nabla g(Y),                                   \label{eqn: admm-ode-y}               \\
        & s^2 \dot{\Lambda} = FX + GY - h.                                      \label{eqn: admm-ode-lambda}  
\end{empheq}    
\end{subequations}
It is important to note that the implicit (Euler) scheme of the system~\eqref{eqn: admm-ode-x} --- \eqref{eqn: admm-ode-lambda}  is exactly the~\texttt{ADMM} algorithm~\eqref{eqn: admm-x} --- \eqref{eqn: admm-lambda}. Furthermore, the system of high-resolution ODEs~\eqref{eqn: admm-ode-x} --- \eqref{eqn: admm-ode-lambda} can capture the effect of the $\lambda$-correction,  which is a significant feature of~\texttt{ADMM}. The small but essential perturbation causes the pair of variable $(X, Y)$ to deviate from the constrained hyperplane $Fx + Gy = h$. We visualize this mathematical observation in~\Cref{fig: admm-low-high}.

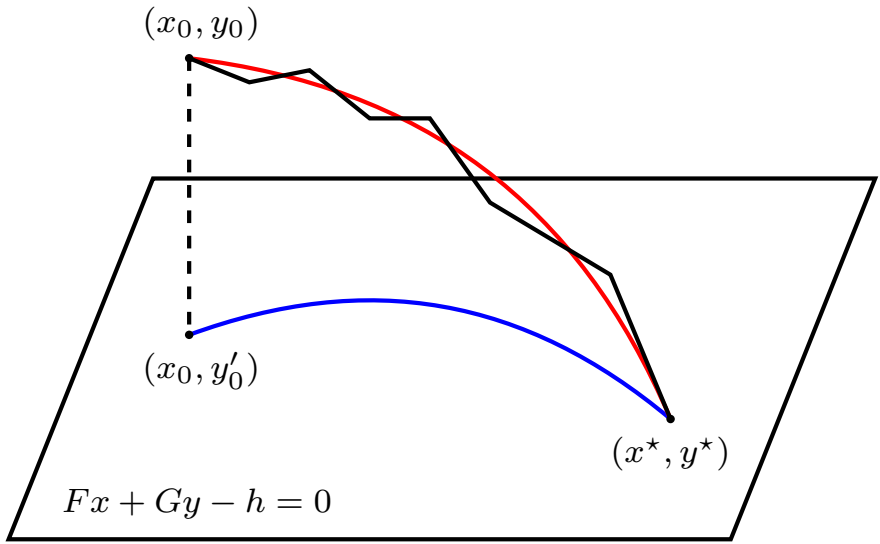
\begin{figure}[htb!]
  \centering
  \scalebox{1.6}{\begin{tikzpicture}[line width=1.0pt].
    \draw[xslant=0.4](-1,0) rectangle (5,3);
    \node at(-0.7,0.3)[right]{\scriptsize\( Fx + Gy - h = 0 \)};
  
    \path[blue] (0.5, 1.7) edge [bend left] (4.5,1);
    \fill(0.5,1.7)circle(1pt)node[below, xshift=0.1cm]{\scriptsize\( (x_0, y_0') \)};

    \path[red] (0.5, 4) edge [bend left] (4.5,1);

    \draw[black] (0.5,4) -- (1.0, 3.8) -- (1.5,3.9)-- (2.0,3.5) -- (2.5,3.5) -- (3.0, 2.8) -- (3.5,2.5) -- (4,2.2) -- (4.5,1);
    \fill(0.5,4)circle(1pt)node[above, xshift=0.1cm]{\scriptsize\( (x_0, y_0) \)};

    \fill(4.5,1)circle(1pt)node[below]{\scriptsize\( (x^\star, y^\star) \)};

    \draw[black,dashed](0.5,4) -- (0.5,1.7); 
  \end{tikzpicture}}
\caption{A schematic diagram of the trajectory of~\texttt{ADMM} with any initial $(x_0, y_0)$ (Black), the continuous limit ODE with any initial $x_0$ (Blue) and the system of high-resolution ODEs with any initial $(x_0, y_0)$ (Red), where the initial $y_0'$ is not arbitrary and required to satisfies $Fx_0 +Gy_0' = h$.}
\label{fig: admm-low-high}
\end{figure}

\item[(2)] Starting from the system of high-resolution ODEs~\eqref{eqn: admm-ode-x} --- \eqref{eqn: admm-ode-lambda},  we can provide a straightforward proof for a quadratic Lyapunov function that decreases monotonically in the continuous case. This calculation can be directly transformed from the continuous system of high-resolution ODEs~\eqref{eqn: admm-ode-x} --- \eqref{eqn: admm-ode-lambda} to the discrete~\texttt{ADMM} algorithm~\eqref{eqn: admm-x} --- \eqref{eqn: admm-lambda}. It is worth noting that the numerical error resulting from the implicit discretization leads to the convergence of~\texttt{ADMM} with a rate of $O(1/N)$.  Furthermore, it has been observed that the numerical error resulting from the implicit discretization still ensures that another discrete Lyapunov function decreases monotonically. Compared to previous results shown in~\citep{he2012convergence, he2015non}, our proofs are principled, succinct, and easy to understand.

\item[(3)] In addition, we further discover that if the objective function $f$ is strongly convex, the iterative sequence of~\texttt{ADMM} converges strongly in terms of average with the following rate
\[
\|\overline{x_N} - x^{\star}\|^2 \leq O\left( \frac{1}{N+1}\right),
\]
where $\overline{x_N}$ takes the form $1/(N+1) \sum_{k=0}^{N} x_k$. This assumption of strong convexity is applicable in practical scenarios, particularly for the generalized~\texttt{Lasso}.
\end{itemize}

%
%

\section{Preliminaries}
\label{sec: prelim}

In this section, we provide a brief description of basic definitions and classical theorems in convex optimization. These definitions and theorems are primarily referenced from classical books, such as~\citep{rockafellar1970convex, rockafellar2009variational} and~\citep{nesterov1998introductory}, and will serve as references for proofs in the subsequent sections. 

At the beginning, let us start with the definition of a convex function and its subgradients and subdifferential. 
\begin{defn}
\label{defn: convex} 
A function $f$ is said to be convex if, for  any $x, y \in \mathbb{R}^d$ and any $\alpha \in [0,1]$, the following inequlaity holds:
\[
f\left( \alpha x + (1 - \alpha)y \right) \leq \alpha f(x) + (1 - \alpha)f(y).
\] 
\end{defn}

\begin{defn}
\label{defn: subgradient}
 A vector $v$ is called a subgradient of $f$ at $x$ if for any $y \in \mathbb{R}^d$, the following inequlity holds:
\[
f(y) \geq f(x) + \langle v, y - x\rangle.
\] 
Additionally, the set of all subgradients of $f$ at $x$ is called the subdifferential of $f$ at $x$, and we note it as $\partial f(x)$. 
\end{defn}

Let $\mathcal{C}(\mathbb{R}^d)$ be the class of continuous functions on $\mathbb{R}^d$. Then, we can denote $\mathcal{F}^0(\mathbb{R}^d) \subset \mathcal{C}(\mathbb{R}^d)$ as the class of continuous and convex functions, and $\mathcal{F}^1(\mathbb{R}^d) \subset \mathcal{F}^0(\mathbb{R}^d)$ as its subclass, the class of differentiable and convex functions. Since the function value is finite everywhere for any $f \in \mathcal{F}^0(\mathbb{R}^d)$, we can conclude that subgradients always exist according to Proposition 8.21 in~\citep{rockafellar1970convex}. We formally state this as the following theorem. 

\begin{theorem}
\label{thm: subgradient-exist}
Let $f \in \mathcal{F}^{0}(\mathbb{R}^d)$. Then for any $x \in \mathbb{R}^{d}$, the subdifferential $\partial f(x)$ is nonempty. 
\end{theorem}

Furthermore, we can provide a sufficient and necessary condition for a minimizer of any continuous and convex function. 

\begin{theorem}
\label{thm: minimizer}
Let $f\in \mathcal{F}^0(\mathbb{R}^d)$. A point $x^{\star}$ is a minimizer of $f$ if and only if $0 \in \partial f(x^{\star})$. 
\end{theorem}

Next, we formally describe two basic properties of subdifferentials and subgradients, which are collected from~\citep{rockafellar1970convex}. 

\begin{theorem}[Theorem 23.8 in \citep{rockafellar1970convex}]
\label{thm: sum-diff}
Let $f_1, f_2 \in\mathcal{F}^0(\mathbb{R}^d)$. Then. for any $x \in \mathbb{R}^d$, the subdifferentials satisfy
\[
\partial(f_1 + f_2)(x) = \partial f_1(x) + \partial f_2(x). 
\]
\end{theorem}

\begin{theorem}[Theorem 25.1 in \citep{rockafellar1970convex}]
\label{thm: subgrad-unique}
Let $f \in \mathcal{F}^1(\mathbb{R}^d)$. Then, the gradient $\nabla f(x)$ is the unique subgradient of $f$ at $x$, so that in particular 
\[
f(y) \geq f(x) + \langle \nabla f(x), y - x\rangle,
\] 
for any $y \in \mathbb{R}^d$.
\end{theorem}

For smooth functions, there still exists some precise characterizations, which we state as the following definition. 
\begin{defn}
\label{defn: L-smooth-strongly}
Let $f \in \mathcal{F}^1(\mathbb{R}^d)$. For any $x, y \in \mathbb{R}^d$, 
\begin{itemize}
\item[(1)] if the gradient of the function $f$ satisfies
\[
\|\nabla f(x) - \nabla f(y)\| \leq L \|x - y\|,
\]
we said that $f$ is $L$-smooth;
\item[(2)] if the function $f$ satisfies 
\[
f(y) \geq f(x) + \langle \nabla f(y), y - x \rangle + \frac{\mu}{2} \|y - x\|^2, 
\]
we said that $f$ is $\mu$-strongly convex. 
\end{itemize}
We note $\mathcal{S}_{\mu,L}^1(\mathbb{R}^d) \subset \mathcal{F}^1(\mathbb{R}^d)$ as its subclass, the class of $L$-smooth and $\mu$-strongly convex functions.
\end{defn}

Finally, let us briefly introduce some basic facts regarding the constrained optimization problem~\eqref{eqn: admm-problem} and~\texttt{ADMM}.  By incorporating the constraint into the objective function, we can obtain the unaugmented Lagrangian for the constrained optimization problem~\eqref{eqn: admm-problem} as 
\begin{equation}
\label{eqn: lagrangian-admm-problem}
L(x,y; \lambda) = f(x) + g(y) + \langle \lambda,  Fx + Gy - h \rangle.
\end{equation}
Based on the dual theorem in convex optimization, we can derive the basic statement expressed in the following theorem. 
\begin{theorem}[Theorem 11.50 in \citep{rockafellar2009variational}]
\label{thm: saddle-equivalence}
The constrained optimization problem~\eqref{eqn: admm-problem} is equivalent to the following unconstrained minimax problem as
\begin{equation}
\label{eqn: minimax}
 \max_{\lambda \in \mathbb{R}^{m}} \min_{ \substack{x \in \mathbb{R}^{d_1} \\y \in \mathbb{R}^{d_2} } } L(x,y; \lambda) = \min_{ \substack{x \in \mathbb{R}^{d_1} \\y \in \mathbb{R}^{d_2} } }\max_{\lambda \in \mathbb{R}^{m}}  L(x,y; \lambda).
\end{equation}
Let $(x^{\star},y^{\star}; \lambda^{\star})$ be the saddle point of the unconstrained minimax problem~\eqref{eqn: minimax}. Furthermore, the unconstrained minimax problem~\eqref{eqn: minimax} is equivalent to the following inequality as
\begin{equation}
\label{eqn: critical-inequality}
 f(x^{\star}) - f(x) + g(y^{\star}) - g(y) \leq \big\langle \lambda^{\star}, F(x - x^{\star}) + G(y - y^{\star})\big\rangle. 
\end{equation}
\end{theorem}
Since the saddle point $(x^{\star},y^{\star}; \lambda^{\star})$ satisfies $Fx^{\star} + Gy^{\star} = h$, it is straightforward for us to derive the equivalence between the minimax problem~\eqref{eqn: minimax} and the inequality~\eqref{eqn: critical-inequality}, which can be established through the bridge inequality $L(x^\star, y^\star; \lambda) \leq L(x^\star, y^\star; \lambda^{\star}) \leq L(x, y; \lambda^{\star})$. In the~\texttt{ADMM} algorithm, the implicit solution in the second iteration~\eqref{eqn: admm-y} can always be obtained directly since the convex function $g$ is assumed to be the $\ell_1$-norm. However, in some cases, such as $f$ being $L$-smooth, the implicit solution in the first iteration~\eqref{eqn: admm-x} may not be obtained directly. Therefore, we still need to assume the class of the objective function $f$ that the implicit solution in the first iteration~\eqref{eqn: admm-x} can be directly obtained and denote it as $\mathcal{R}(\mathbb{R}^d)$. 

%
%



%

\section{Perspective from the system of high-resolution ODEs }
\label{sec: ode}

In this section, we first derive the system of high-resolution ODEs~\eqref{eqn: admm-ode-x} --- \eqref{eqn: admm-ode-lambda}  for the~\texttt{ADMM} algorithm by use of dimensional analysis. Then, we employ Lyapunov analysis to show the convergence behaviors of the system of high-resolution ODEs. 

\subsection{Derivation of the system of high-resolution ODEs}
\label{subsec: derivation-ode}

To derive the system of high-resolution ODEs, it is necessary to assume that both the objective functions and the solutions are sufficiently smooth. Furthermore, we set the parameter $s$ as the step size, which serves as a bridge between the~\texttt{ADMM} algorithm and the continuous system. Let $t_{k}=ks$, $(k=0,1,2, \ldots)$ and take the ansatz $x_k = X(t_k), \; y_k = Y(t_k), \; \lambda_k = \Lambda(t_k)$. Performing a Taylor expansion in powers of $s$, we can obtain 
\begin{equation}
\label{eqn: taylor}
\left\{ \begin{aligned}
        & x_{k} = X(t_{k}) = X(t_{k+1}) - s\dot{X}(t_{k+1})  + O\left(s^2\right) = x_{k+1} - s\dot{X}(t_{k+1})  + O\left(s^2\right) \\
        & y_{k} = Y(t_{k}) \;= Y(t_{k+1}) - s\dot{Y}(t_{k+1}) + O\left(s^2\right) \; = y_{k+1} - s\dot{Y}(t_{k+1})  + O\left(s^2\right)    \\
        & \lambda_{k} = \Lambda(t_{k}) \;= \Lambda(t_{k+1}) \;- s\dot{\Lambda}(t_{k+1})  + O\left(s^2\right) \; = \lambda_{k+1} - s\dot{\Lambda}(t_{k+1})  + O\left(s^2\right). 
        \end{aligned} \right.
\end{equation}

For any $f \in \mathcal{F}^1(\mathbb{R}^d) \cap \mathcal{R}(\mathbb{R}^d)$ and  $g \in \mathcal{F}^1(\mathbb{R}^d)$,  the first two iterations of~\texttt{ADMM},~\eqref{eqn: admm-x} and~\eqref{eqn: admm-y}, can be explicitly written as:
\begin{subequations}
\label{eqn: admm-explicit}
\begin{empheq}[left=\empheqlbrace]{align}
        & s\nabla f(x_{k+1}) +  F^{\top}(Fx_{k+1} + G y_k - h + s \lambda_k) = 0,              \label{eqn: admm-explicit-x} \\
        & s\nabla g(y_{k+1}) +  G^{\top}(F x_{k+1} + G y_{k+1} - h + s \lambda_k) = 0,     \label{eqn: admm-explicit-y} 
\end{empheq}    
\end{subequations} 
By utilizing the third iteration~\eqref{eqn: admm-lambda}, we can reformulate the first iteration~\eqref{eqn: admm-explicit-x} as
\begin{equation}
\label{eqn: admm-explicit-new-x}
F^{\top}G\cdot \frac{y_{k+1} - y_{k}}{s}    =  \nabla f(x_{k+1}) + F^{\top} \lambda_{k+1}    
\end{equation}
Substituting~\eqref{eqn: taylor} into~\eqref{eqn: admm-explicit-new-x}, we have
\[
F^{\top}G \dot{Y}(t_{k+1})  - O(s) =  \nabla f(X(t_{k+1})) + F^{\top} \Lambda(t_{k+1} ),  
\]
where we can derive the ODE~\eqref{eqn: admm-ode-x} by considering the $O(1)$ terms and ignoring the $O(s)$ terms.  If we further consider the $O(s)$ terms and ignore the high-order terms, it corresponds to an overdamped system.  In an overdamped system, any inertial effects are damped out quickly by the effects of viscous, frictional, or other damping forces, and can be neglected~\citep{adams2013large}. Therefore, the ODE~\eqref{eqn: admm-ode-x} is a reasonable approximation for the first iteration~\eqref{eqn: admm-explicit-new-x}. Furthermore, the first iteration~\eqref{eqn: admm-explicit-new-x} is exactly the implicit discretization of the ODE~\eqref{eqn: admm-ode-x}. 

By utilizing the third iteration~\eqref{eqn: admm-lambda}, we can rephrase the second iteration~\eqref{eqn: admm-explicit-x} as an algebraic equation:
\begin{equation}
\label{eqn: admm-explicit-new-y}  
\nabla g(y_{k+1}) + G^{\top} \lambda_{k+1} =  0,   
\end{equation}
which corresponds to the algebraic equation~\eqref{eqn: admm-ode-y}. Finally, substituting~\eqref{eqn: taylor} into the third iteration~\eqref{eqn: admm-lambda}, we have
\[
s^2\dot{\Lambda}(t_{k+1})  -  O\left(s^3\right) = FX(t_{k+1}) + GY(t_{k+1}) - h
\]
where the high-resolution ODE~\eqref{eqn: admm-ode-lambda} can be derived by including $O(s^2)$ terms. However, if we only include $O(s)$ terms and ignore $O(s^2)$ terms, the high-resolution ODE~\eqref{eqn: admm-ode-lambda} retrogrades to the algebraic equation of the constrained hyperplane $FX + GY - h = 0$. Furthermore, the system of high-resolution ODEs~\eqref{eqn: admm-ode-x} --- \eqref{eqn: admm-ode-lambda} reduces to the continuous-limit ODE~\eqref{eqn: ode-low}, so it is also referred to as the low-resolution ODE of~\texttt{ADMM}.


%
%
%

\subsection{Convergence of the system of high-resolution ODEs}
\label{subsec: convergence-ode}

Let $(X,Y; \Lambda)$ be the solution to the system of high-resolution ODEs~\eqref{eqn: admm-ode-x} --- \eqref{eqn: admm-ode-lambda} and $(x,y; \lambda)$ be any point in $\mathbb{R}^{d_1} \times \mathbb{R}^{d_2} \times \mathbb{R}^{m}$. To analyze the convergence of the system of high-resolution ODEs, we construct the following Lyapunov function as
\begin{equation}
\label{eqn: lyapunov-ode}
\mathcal{E}(t) =  \frac{1}{2s} \left\|G(Y - y)\right\|^2 + \frac{s}{2} \left\| \Lambda - \lambda \right\|^2.
\end{equation}

\begin{lemma}
\label{lem: ode}
Let $f \in \mathcal{F}^1(\mathbb{R}^d) \cap \mathcal{R}(\mathbb{R}^d)$ and  $g \in \mathcal{F}^1(\mathbb{R}^d)$. For any step size $s$, and with $(X,Y; \Lambda) = (X(t), Y(t); \Lambda(t))$ being the solution to the system of high-resolution ODEs~\eqref{eqn: admm-ode-x} --- \eqref{eqn: admm-ode-lambda}, the Lyapunov function~\eqref{eqn: lyapunov-ode} satisfies the following inequality as
\begin{align}
\frac{d \mathcal{E}}{dt} \leq \frac{1}{s} \bigg( f(x) - f(X) + g(y) - g(Y) &+ \big\langle \Lambda - G\dot{Y}, F(x-x^{\star}) + G(y - y^{\star}) \big\rangle \nonumber \\
 & - \big\langle \lambda, F (X - x^{\star}) + G (Y-y^{\star}) \big\rangle \bigg). \label{eqn: derivative-inequality}
\end{align}
\end{lemma}

\begin{proof}[Proof of~\Cref{lem: ode}]
With the high-resolution ODE~\eqref{eqn: admm-ode-lambda}, we can calculate the time derivative of the Lyapunov function~\eqref{eqn: lyapunov-ode} as 
\begin{align}
\frac{d \mathcal{E}}{dt} & = \frac{1}{s} \big\langle G\dot{Y}, G(Y - y) \big\rangle + s \big\langle \dot{\Lambda}, \Lambda - \lambda \big\rangle                \nonumber \\
                                     & = \frac{1}{s} \left( \big\langle G\dot{Y}, G(Y - y) \big\rangle + \big\langle  F (X - x^{\star}) + G (Y-y^{\star}), \Lambda - \lambda \big\rangle \right). \label{eqn: lem-ode-1}
\end{align}
Since $f \in \mathcal{F}^1(\mathbb{R}^d)$ and $g \in \mathcal{F}^1(\mathbb{R}^d)$, we can derive the two convex inequalities from~\Cref{thm: subgrad-unique} as
\begin{align}
& f(x) - f(X) \geq \big\langle \nabla f(X), x - X \big\rangle = \big\langle-F^{\top}\Lambda +  F^{\top}G\dot{Y}, x- X \big\rangle, \label{eqn: lem-ode-2} \\
& g(y) - g(Y)\; \geq \big\langle \nabla g(Y), y - Y \big\rangle \; = \big\langle  -G^{\top}\Lambda, y - Y \big\rangle, \label{eqn: lem-ode-3}
\end{align}
where the two equalities directly follows the first two iterations of~\texttt{ADMM},~\eqref{eqn: admm-ode-x} and \eqref{eqn: admm-ode-y}. Combining~\eqref{eqn: lem-ode-2} and~\eqref{eqn: lem-ode-3}, we further have
\begin{equation}
\label{eqn: lem-ode-4}
f(x) - f(X) + g(y) - g(Y) + \big\langle G\dot{Y}, F(X - x) \big\rangle - \big\langle \Lambda, F(X - x) + G(Y - y) \big\rangle \geq 0.
\end{equation}
Plugging~\eqref{eqn: lem-ode-4} into the derivative inequality~\eqref{eqn: lem-ode-1}, we have
\begin{align*}
\frac{d \mathcal{E}}{dt}             & \leq \frac{1}{s} \bigg( f(x) - f(X) + g(y) - g(Y) + \big\langle \Lambda, F(x-x^{\star}) + G(y - y^{\star}) \big\rangle  \\
                                     &\mathrel{\phantom{\leq \frac{1}{s}()}} + \;\big\langle G\dot{Y}, F(X-x) + G(Y-y) \big\rangle  -  \big\langle \lambda, F (X - x^{\star}) + G (Y-y^{\star}) \big\rangle \bigg) \\
                                     & =  \frac{1}{s} \bigg( f(x) - f(X) + g(y) - g(Y) + \big\langle \Lambda - G\dot{Y}, F(x-x^{\star}) + G(y - y^{\star}) \big\rangle  \\
                                     &\mathrel{\phantom{\leq \frac{1}{s}()}}  - \; \big\langle \lambda, F (X - x^{\star}) + G (Y-y^{\star}) \big\rangle  \bigg) + \big\langle G\dot{Y},\dot{\Lambda} \big\rangle  \\
                                     & =  \frac{1}{s} \bigg( f(x) - f(X) + g(y) - g(Y) + \big\langle \Lambda - sG\dot{Y}, F(x-x^{\star}) + G(y - y^{\star}) \big\rangle  \\
                                     &\mathrel{\phantom{\leq \frac{1}{s}(}}  - \; \big\langle \lambda, F (X - x^{\star}) + G (Y-y^{\star}) \big\rangle \bigg) - s\dot{Y}^{\top}\nabla^2 g(Y) \dot{Y}, 
\end{align*} 
where the two equalities directly follow the high-resolution ODEs,~\eqref{eqn: admm-ode-y} and~\eqref{eqn: admm-ode-lambda}. Since the inequality $\dot{Y}^{\top}\nabla^2 g(Y) \dot{Y} \geq 0$ always holds, we complete the proof. 
\end{proof}

In the given context, the time average of a variable, denoted as $X$, within the time interval $[0,t]$, is defined as 
\[
\overline{X} = \frac{1}{t}\int_{0}^{t} X(s)ds. 
\]
It is important to mention that the parameter $\lambda$ in the derivative inequality~\eqref{eqn: derivative-inequality} is arbitrary. Therefore, we have the flexibility to set it to zero, which directly determines the convergence rate of the time average in the weak sense by utilizing~\Cref{lem: ode}. We state it rigorously as follows.  
\begin{theorem}
\label{thm: ode-average}
Under the assumptions of~\Cref{lem: ode}, the time average $(\overline{X}, \overline{Y},\overline{\Lambda - G\dot{Y}})$ converges weakly to $(x^{\star},y^{\star};\lambda^{\star})$  with the following rate as
\begin{equation}
\label{eqn: ode-average}
f(\overline{X}) - f(x) + g(\overline{Y}) - g(y) + \big\langle \overline{\Lambda - G\dot{Y}}, F(x-x^{\star}) + G(y - y^{\star}) \big\rangle \leq \frac{\|G(y_0 - y)\|^2 + s^2 \|\lambda_0\|^2}{2t},
\end{equation}
for any initial $(x_0, y_0; \lambda_0) \in \mathbb{R}^{d_1} \times \mathbb{R}^{d_2} \times \mathbb{R}^{m}$. 
\end{theorem}

When the point $(x,y; \lambda)$ is set as the saddle point $(x^{\star}, y^{\star}; \lambda^{\star})$, the Lyapunov function~\eqref{eqn: lyapunov-ode} takes the following form as
\begin{equation}
\label{eqn: lypunov-ode-star}
\mathcal{E}(t) =  \frac{1}{2s} \left\|G(Y - y^{\star})\right\|^2 + \frac{s}{2} \left\| \Lambda - \lambda^{\star} \right\|^2.
\end{equation}
Plugging the saddle point $(x^{\star}, y^{\star}; \lambda^{\star})$ into the derivative inequality~\eqref{eqn: derivative-inequality}, we can obtain 
\begin{equation}
\label{eqn: ode-convergence1}
\frac{d \mathcal{E}}{dt} \leq \frac{1}{s} \bigg( f(x^{\star}) - f(X) + g(y^{\star}) - g(Y)  - \big\langle \lambda^{\star}, F (X - x^{\star}) + G (Y-y^{\star}) \big\rangle \bigg). 
\end{equation}
With the help of~\Cref{thm: saddle-equivalence}, we can establish that the Lyapunov function decreases monotonally. The rigorous representation is stated as follows.
\begin{theorem}
\label{thm: ode-convegence}
Under the assumption of~\Cref{lem: ode}, the Lyapunov function given by~\eqref{eqn: lypunov-ode-star} decreases monotonically. 
\end{theorem}

In addition, considering the generalized~\texttt{Lasso},  the objective function $f$ may be $\mu$-strongly convex. Assuming that the objective function $f$ does indeed possess $\mu$-strongly convexity, we can deduce that the time average of $X$ over the time interval $[0,t]$ converges strongly with a rate of $O(1/t)$. The rigorous representation is stated as follows. 

\begin{theorem}
\label{thm: ode-strongly}
Under the assumptions of~\Cref{lem: ode}, if we further assume the objective function satisfies $f \in \mathcal{S}^1_{\mu,L}(\mathbb{R}^d)$, the time average $\overline{X}$ converges with the following rate as
\[
\left\|\overline{X} - x^{\star}\right\|^2 \leq \frac{\|x_0 - x^{\star}\|^2 + s^2\|\lambda_0 - \lambda^{\star}\|^2}{\mu t},
\]
for any initial $(x_0, y_0; \lambda_0) \in \mathbb{R}^{d_1} \times \mathbb{R}^{d_2} \times \mathbb{R}^{m}$. 
\end{theorem}

\begin{proof}[Proof of~\Cref{thm: ode-strongly}]
When the objective function is assumed to be $f \in \mathcal{S}^1_{\mu,L}(\mathbb{R}^d)$, we can tighten the convex inequality~\eqref{eqn: lem-ode-2} as  
\begin{equation}
\label{eqn: ode-strong-convex}
f(x^{\star}) - f(X) \geq \big\langle-F^{\top}\Lambda +  F^{\top}G\dot{Y}, x^{\star} - X \big\rangle + \frac{\mu}{2} \left\|x^{\star} - X\right\|^2. 
\end{equation}
By utilizing the inequality~\eqref{eqn: ode-strong-convex}, and following the same process used to obtain~\Cref{thm: ode-convegence}, we can deduce that the time derivative of the Lyapunov function~\eqref{eqn: lyapunov-ode} satisfies the following inequality:
\[
\frac{d \mathcal{E}}{dt} \leq - \frac{\mu}{2s}  \left\|X - x^{\star}\right\|^2.
\]
Hence, by taking the time average of $X$, we complete the proof with the convexity of the $\ell_2$-norm square.   
\end{proof}

%
%
%
%
%
%
%
%
%
%
%
%
%
%
%

\section{Convergence rates of ADMM}
\label{sec: admm}

In~\Cref{subsec: convergence-ode}, we have delved into the convergence behavior of the system of high-resolution ODEs~\eqref{eqn: admm-ode-x} --- \eqref{eqn: admm-ode-lambda}. Now, we extend our exploration straightforward from the continuous perspective and shift our focus to the discrete~\texttt{ADMM} algorithm~\eqref{eqn: admm-x} --- \eqref{eqn: admm-lambda}. In this context, we can discover the numerical phenomenon depicted in~\Cref{subsec: error-implicit},  where the numerical error arising from the implicit discretization significantly affects the convergence rate. 

Let $(x,y; \lambda) \in \mathbb{R}^{d_1} \times \mathbb{R}^{d_2} \times \mathbb{R}^{m}$ be any point and $\{(x_k,y_k;\lambda_k)\}_{k=0}^{\infty}$ be the iterative sequence of~\texttt{ADMM}~\eqref{eqn: admm-x} --- \eqref{eqn: admm-lambda}. We extend the continuous Lyapunov function~\eqref{eqn: lyapunov-ode} to its discrete form as
\begin{equation}
\label{eqn: admm-lyapunov}
\mathcal{E}(k) = \frac{1}{2s} \left\|G(y_k - y)\right\|^2 + \frac{s}{2} \left\| \lambda_k - \lambda \right\|^2,
\end{equation}
which results in the following lemma.

\begin{lemma}
\label{lem: admm}
Let $f \in \mathcal{F}^0(\mathbb{R}^d) \cap \mathcal{R}(\mathbb{R}^d)$ and  $g \in \mathcal{F}^0(\mathbb{R}^d)$. Then, the discrete Lyapunov function~\eqref{eqn: admm-lyapunov} satisfies the following inequality as
\begin{align}
\mathcal{E}&(k+1)  - \mathcal{E}(k) \nonumber \\
& \leq f(x) - f(x_{k+1}) + g(y) - g(y_{k+1}) + \left\langle \lambda_{k+1} - \frac{G(y_{k+1} - y_{k})}{s}, F(x - x^{\star}) + G(y - y^{\star}) \right\rangle \nonumber \\
                                       & \mathrel{\phantom{\leq}} - \big\langle \lambda,  F(x_{k+1} - x^{\star}) + G(y_{k+1} - y^{\star}) \big\rangle - \underbrace{\left( \frac{1}{2s} \left\|G(y_{k+1} - y_{k})\right\|^2 + \frac{s}{2} \left\| \lambda_{k+1} - \lambda_{k}\right\|^2 \right)}_{\mathbf{NE}},  \label{eqn: iterative-inequality}
\end{align}
where $\mathbf{NE}$ represents the numerical error resulting from the implicit scheme. 
\end{lemma}

\begin{proof}[Proof of~\Cref{lem: admm}]
With the discrete Lyapunov function~\eqref{eqn: admm-lyapunov}, we can calculate its iterative difference as
\begin{align}
\mathcal{E}(k+1) - \mathcal{E}(k)    & =  \underbrace{\frac{1}{s} \big\langle G(y_{k+1} - y_{k}), G(y_{k+1} - y) \big\rangle +  s\big\langle \lambda_{k+1} - \lambda_{k}, \lambda_{k+1} - \lambda \big\rangle}_{\mathbf{I}} \nonumber \\
                                                     & \mathrel{\phantom{=}} - \underbrace{ \left( \frac{1}{2s} \left\|G(y_{k+1} - y_{k})\right\|^2 + \frac{s}{2} \left\| \lambda_{k+1} - \lambda_{k}\right\|^2 \right)}_{\mathbf{NE}}, \label{eqn: admm-iterative-diff}
\end{align}
where $\mathbf{I}$ corresponds to the one obtained from the continuous perspective, and $\mathbf{NE}$ represents the numerical error resulting from the implicit discretization. By utilizing the third iteration of~\texttt{ADMM}~\eqref{eqn: admm-lambda}, we can reformulate $\mathbf{I}$ as 
\begin{equation}
\label{eqn: admm-lem-1}
\mathbf{I} = \frac{1}{s} \big\langle G(y_{k+1} - y_{k}), G(y_{k+1} - y) \big\rangle +  \big\langle  F(x_{k+1} - x^{\star}) + G(y_{k+1} - y^{\star}), \lambda_{k+1} - \lambda \big\rangle.
\end{equation}
In the case of the discrete~\texttt{ADMM} algorithm, both $f \in \mathcal{F}^0(\mathbb{R}^d) \cap \mathcal{R}(\mathbb{R}^d)$ and  $g \in \mathcal{F}^0(\mathbb{R}^d)$ are nonsmooth.  As a result, we need to apply~\Cref{thm: minimizer} and~\Cref{thm: sum-diff} to explicitly reformulate the first two iterations of~\texttt{ADMM},~\eqref{eqn: admm-x} and~\eqref{eqn: admm-y}, which can be expressed as
\begin{subequations}
\label{eqn: admm-nonsmooth-explicit}
\begin{empheq}[left=\empheqlbrace]{align}
        & s\partial f(x_{k+1}) +  F^{\top}(Fx_{k+1} + G y_k - h + s \lambda_k) \ni  0,              \label{eqn: admm-nonsmooth-explicit-x} \\
        & s\partial g(y_{k+1}) +  G^{\top}(F x_{k+1} + G y_{k+1} - h + s \lambda_k) \ni  0.         \label{eqn: admm-nonsmooth-explicit-y} 
\end{empheq}    
\end{subequations}
Furthermore, plugging the third iteration of~\texttt{ADMM}~\eqref{eqn: admm-lambda} into both the first two iterations,~\eqref{eqn: admm-nonsmooth-explicit-x} and~\eqref{eqn: admm-nonsmooth-explicit-y}, we have
\begin{subequations}
\label{eqn: admm-nonsmooth-explicit-1}
\begin{empheq}[left=\empheqlbrace]{align}
        & \partial f(x_{k+1}) - F^{\top}G\cdot \frac{y_{k+1} - y_{k}}{s} + F^{\top} \lambda_{k+1} \ni 0 ,    \label{eqn: admm-nonsmooth-explicit-x-1} \\
        & \partial g(y_{k+1}) + G^{\top} \lambda_{k+1} \ni 0.                                                \label{eqn: admm-nonsmooth-explicit-y-1} 
\end{empheq}    
\end{subequations}
With the help of~\eqref{eqn: admm-nonsmooth-explicit-x-1} and~\eqref{eqn: admm-nonsmooth-explicit-y-1}, we can derive two convex inequalities from~\Cref{defn: subgradient} and~\Cref{thm: sum-diff} as
\begin{align}
& f(x) - f(x_{k+1}) \geq  \left\langle F^{\top}G\cdot \frac{y_{k+1} - y_{k}}{s} - F^{\top} \lambda_{k+1}, x - x_{k+1} \right\rangle,           \label{eqn: lem-admm-2} \\
& g(y) - g(y_{k+1}) \;\geq  \big\langle  -G^{\top}\lambda_{k+1}, y - y_{k+1} \big\rangle. \label{eqn: lem-admm-3}
\end{align}
By adding these two convex inequalities~\eqref{eqn: lem-admm-2} and~\eqref{eqn: lem-admm-3} to~\eqref{eqn: admm-lem-1}, we obtain
\begin{align*}
\mathbf{I} \leq & f(x) - f(x_{k+1}) + g(y) - g(y_{k+1}) + \left\langle \lambda_{k+1} - \frac{G(y_{k+1} - y_{k})}{s}, F(x - x^{\star}) + G(y - y^{\star}) \right\rangle  \\
                & - \big\langle \lambda,  F(x_{k+1} - x^{\star}) + G(y_{k+1} - y^{\star}) \big\rangle  + \big\langle y_{k+1} - y_{k}, G^{T}(\lambda_{k+1} - \lambda_{k}) \big\rangle,
\end{align*}
where the last term $\big\langle y_{k+1} - y_{k}, G^{T}(\lambda_{k+1} - \lambda_{k}) \big\rangle \leq 0$ always holds due to the basic convex inequality from~\Cref{defn: subgradient}. Therefore, the proof is complete. 
\end{proof}

In the given context, we denote the average of an iterative sequence. Taking $\{x_{k}\}_{k=0}^{\infty}$ for example, the iterative average is formulated as 
\[
\overline{x_{N}} = \frac{1}{N+1} \sum_{k=0}^{N} x_k. 
\]
Since $\lambda$ in the inequality regarding the iterative difference~\eqref{eqn: iterative-inequality} is arbitrary, we can set it to zero. Furthermore, by utilizing~\Cref{lem: admm}, we can determine the rate that the average sequence $\{ \overline{x_N} \}_{N=0}^{\infty}$ converges in the weak sense. We state this result rigorously as follows.  
\begin{theorem}
\label{thm: admm-average}
Under the assumption of~\Cref{lem: admm}, the iterative average $(\overline{x_N}, \overline{x_N},\overline{\lambda_N} - G(\overline{y_{N+1} - y_{N}})/s)$ converges weakly to $(x^{\star},y^{\star};\lambda^{\star})$  with the following rate as
\begin{multline}
\label{eqn: admm-average}
f(\overline{x_N}) - f(x) + g(\overline{y_{N}}) - g(y) \\ + \left\langle \overline{\lambda_N} - \frac{G(\overline{y_{N+1} - y_{N}})}{s}, F(x-x^{\star}) + G(y - y^{\star}) \right\rangle \leq \frac{\|G(y_0 - y)\|^2 + s^2 \|\lambda_0\|^2}{2s(N+1)},
\end{multline}
for any initial $(x_0, y_0; \lambda_0) \in \mathbb{R}^{d_1} \times \mathbb{R}^{d_2} \times \mathbb{R}^{m}$. 
\end{theorem}

By substituting the saddle point $(x^{\star}, y^{\star}; \lambda^{\star})$ into the Lyapunov function~\eqref{eqn: admm-lyapunov}, we deduce it to have the following form as
\begin{equation}
\label{eqn: lypunov-admm-star}
\mathcal{E}(k) =  \frac{1}{2s} \left\|G(y_k - y^{\star})\right\|^2 + \frac{s}{2} \left\| \lambda_k - \lambda^{\star} \right\|^2.
\end{equation}
Furthermore, the iterative difference~\eqref{eqn: iterative-inequality} can be reduced to satisfy the following inequality as 
\begin{align}
\mathcal{E}(k+1) - \mathcal{E}(k) \leq &  f(x^{\star}) - f(x_{k+1}) + g(y^{\star}) - g(y_{k+1})  - \big\langle \lambda^{\star}, F (x_{k+1} - x^{\star}) + G ( y_{k+1} - y^{\star}) \big\rangle        \nonumber \\ 
                                       & - \underbrace{ \left( \frac{1}{2s} \left\|G(y_{k+1} - y_{k})\right\|^2 + \frac{s}{2} \left\| \lambda_{k+1} - \lambda_{k}\right\|^2 \right)}_{\mathbf{NE}}  \label{eqn: admm-convergence1}
\end{align}
With the help of~\Cref{thm: saddle-equivalence}, we can observe from the inequality~\eqref{eqn: admm-convergence1} that the iterative difference can be negatively bounded by the numerical error resulting from the implicit scheme, while the time derivative is only no more than zero for the continuous case. Therefore, we can derive the convergence rates for the discrete~\texttt{ADMM} algorithm. The rigorous representation is stated as follows. 

%
\begin{theorem}
\label{thm: admm-convegence-rate}
Under the assumption of~\Cref{lem: admm}, the iterative sequence $\{(x_k,y_k; \lambda_k)\}_{k=0}^{\infty}$ generated by~\texttt{ADMM}~\eqref{eqn: admm-x} --- \eqref{eqn: admm-lambda} converges to $(x^{\star},y^{\star};\lambda^{\star})$ with the following rate as
\begin{equation}
\label{eqn:admm-convegence-rate}
\left\{ \begin{aligned}
        & \frac{1}{N+1}\sum_{k=0}^{N} \left( \left\|G(y_{k+1} - y_{k})\right\|^2 + s^2 \left\| \lambda_{k+1} - \lambda_{k} \right\|^2 \right) \leq \frac{\left\|G(y_0 - y^{\star})\right\|^2 + s^2 \left\| \lambda_0 - \lambda^{\star} \right\|^2}{N+1}, \\
        & \min_{0 \leq k \leq N}  \left( \left\|G(y_{k+1} - y_{k})\right\|^2 + s^2 \left\| \lambda_{k+1} - \lambda_{k} \right\|^2 \right) \leq \frac{\left\|G(y_0 - y^{\star})\right\|^2 + s^2 \left\| \lambda_0 - \lambda^{\star} \right\|^2}{N+1},
        \end{aligned} \right.
\end{equation}
for any initial $(x_0, y_0; \lambda_0) \in \mathbb{R}^{d_1} \times \mathbb{R}^{d_2} \times \mathbb{R}^{m}$. 
\end{theorem}

Similarly, if we further make the assumption that the objective function satisfies $f \in \mathcal{S}^1_{\mu,L}(\mathbb{R}^d)$, the convex inequality~\eqref{eqn: lem-admm-2} can become tighter, resulting in the inequality
\[
f(x) - f(x_{k+1}) \geq  \left\langle F^{\top}G\cdot \frac{y_{k+1} - y_{k}}{s} - F^{\top} \lambda_{k+1}, x - x_{k+1} \right\rangle + \frac{\mu}{2} \|x - x_{k+1}\|^2.
\]
Following the same process used to obtain~\Cref{thm: admm-convegence-rate}, we can derive that the iterative difference satisfies the following inequality as
\[
\mathcal{E}(k+1) - \mathcal{E}(k) \leq - \frac{\mu}{2}  \left\|x_{k+1} - x^{\star}\right\|^2 - \underbrace{\left( \frac{1}{2s} \left\|G(y_{k+1} - y_{k})\right\|^2 + \frac{s}{2} \left\| \lambda_{k+1} - \lambda_{k}\right\|^2 \right)}_{\mathbf{NE}}.
\]
Combining with the convexity of the $\ell_2$-norm square, we can derive the convergence rate of the iterative average sequence as stated in the following theorem. 
\begin{theorem}
\label{thm: admm-strongly}
Under the assumption of~\Cref{lem: admm}, if we further assume the objective function $f \in \mathcal{S}^1_{\mu,L}(\mathbb{R}^d)$, the iterative average $\{\overline{x_N}\}_{N=0}^{\infty}$ generated by~\texttt{ADMM}~\eqref{eqn: admm-x} --- \eqref{eqn: admm-lambda} converges to $(x^{\star},y^{\star};\lambda^{\star})$ with the following rate as
\begin{equation}
\label{eqn:admm-convegence-rate-strongly}
\left\|\overline{x_{N+1}} - x^{\star}\right\|^2 \leq \frac{\|x_0 - x^{\star}\|^2 + s^2\|\lambda_0 - \lambda^{\star}\|^2}{\mu s (N+1)},
\end{equation}
for any initial $(x_0, y_0; \lambda_0)$. 
\end{theorem}

\section{Monotonicity}
\label{sec: monotonicity}

In this section, our aim is to investigate the monotonicity of the numerical error, as shown in~\eqref{eqn: iterative-inequality}. The numerical error is given by $\mathbf{NE} = (1/2s)\left\|G(y_{k+1} - y_{k})\right\|^2 + (s/2)) \left\| \lambda_{k+1} - \lambda_{k} \right\|^2$, and we are specifically interested in observing how it changes as the iteration $k$ increases. If the numerical error $\mathbf{NE}$ decreases monotonically as increasing $k$, it indicates that the convergence rate mentioned as $O(1/N)$  in~\eqref{eqn:admm-convegence-rate} can be further enhanced towards the last iterate. This improvement is particularly significant when compared to the convergence rates achieved in terms of average and minimization.

In a similar manner, we will begin by considering the system of high-resolution ODEs~\eqref{eqn: admm-ode-x} --- \eqref{eqn: admm-ode-lambda}. When dealing smooth functions, where $f \in \mathcal{F}^2(\mathbb{R}^d) \cap \mathcal{R}(\mathbb{R}^d)$ and  $g \in \mathcal{F}^2(\mathbb{R}^d)$, we can derive the following inequality as
\begin{equation}
\label{eqn: mono-con-1}
\big\langle \dot{X}, \nabla^2 f(X) \dot{X} \big\rangle + \big\langle \dot{Y}, \nabla^2 g(Y) \dot{Y} \big\rangle \geq 0.
\end{equation}
Taking the time derivative for each ODE in the system~\eqref{eqn: admm-ode-x} --- \eqref{eqn: admm-ode-lambda}, the convex inequality~\eqref{eqn: mono-con-1} can be transformed as 
\begin{equation}
\label{eqn: mono-con-2}
\big\langle F\dot{X}, G\ddot{Y} \big\rangle \geq s^{2} \big\langle \ddot{\Lambda},  \dot{\Lambda} \big\rangle.
\end{equation}
The continuous version of the numerical error $\mathbf{NE}$, as mentioned in~\eqref{eqn: iterative-inequality}, is considered as a Lyapunov function denoted by 
\begin{equation}
\label{eqn: lyapunov-monotone}
\mathcal{E}(t) = \frac{s}{2}\|G\dot{Y}\|^2 + \frac{s^3}{2}\|\dot{\Lambda}\|^2.
\end{equation}
By utilizing the high-resolution ODE~\eqref{eqn: admm-ode-lambda} and the inequality~\eqref{eqn: mono-con-2}, we can deduce that the time derivative satisfies 
\begin{equation}
\label{eqn: mono-con-3}
\frac{d\mathcal{E}}{dt} = s\big\langle G\ddot{Y}, G\dot{Y} \big\rangle + s^3 \big\langle \dot{\Lambda}, \ddot{\Lambda} \big\rangle \leq s^3\big\langle G\ddot{Y}, \ddot{\Lambda} \big\rangle.
\end{equation}
However, since the sign of $s^3\big\langle G\ddot{Y}, \ddot{\Lambda} \big\rangle$ is unknown, we cannot determine the monotonicity of the continuous Lyapunov function~\eqref{eqn: lyapunov-monotone} solely based on the derivative inequality~\eqref{eqn: mono-con-3}. Nevertheless, when we employ the implicit scheme for the system of high-resolution ODEs~\eqref{eqn: admm-ode-x} --- \eqref{eqn: admm-ode-lambda}, there is a possibility that the numerical error, as discussed in~\Cref{subsec: error-implicit} and elaborated in~\Cref{sec: admm}, can offset the term $s^3\big\langle G\ddot{Y}, \ddot{\Lambda} \big\rangle$. Upon closer examination of the Lyapunov function~\eqref{eqn: lyapunov-monotone}, we can observe that the numerical error approximates $- (s^2/2)\|G\dot{Y}\|^2 - s^4/2\|\dot{\Lambda}\|^2$, which can act as a perfect square and potentially counterbalance the term $s^3\big\langle G\ddot{Y}, \ddot{\Lambda} \big\rangle$.

Next, we extend the aforementioned intuition from the system of high-resolution ODEs to the discrete~\texttt{ADMM} algorithm~\eqref{eqn: admm-x} --- \eqref{eqn: admm-lambda}. For any $f \in \mathcal{F}^0(\mathbb{R}^d) \cap \mathcal{R}(\mathbb{R}^d)$ and  $g \in \mathcal{F}^0(\mathbb{R}^d)$, we can derive that the iterative sequence $\{(x_k,y_k;\lambda_k)\}_{k=0}^{\infty}$ satisfies the convex inequalities  from~\Cref{defn: subgradient} as
\begin{align}
\big \langle G(y_{k+2} - y_{k+1}) - G(y_{k+1} - y_{k}), F(x_{k+2} - x_{k+1}) \big \rangle &- s \big \langle \lambda_{k+2} - \lambda_{k+1}, F(x_{k+2} - x_{k+1}) \big\rangle \geq 0,   \label{eqn: mono-admm-1} \\
                                                                                      &-s \big \langle \lambda_{k+2} - \lambda_{k+1}, G(y_{k+2} - y_{k+1}) \big \rangle \geq 0.  \label{eqn: mono-admm-2}
\end{align}
Combining~\eqref{eqn: mono-admm-1} and~\eqref{eqn: mono-admm-2} and incorporating the third iteration~\eqref{eqn: admm-lambda}, we can obtain
\begin{multline}
\big \langle G(y_{k+2} - y_{k+1}) - G(y_{k+1} - y_{k}), F(x_{k+2} - x_{k+1}) \big \rangle \\ 
     \geq   s^{2} \big \langle \lambda_{k+2} - \lambda_{k+1}, (\lambda_{k+2} - \lambda_{k+1}) - (\lambda_{k+1} - \lambda_{k}) \big\rangle. \label{eqn: mono-admm-3}
\end{multline}
As indicated in~\eqref{eqn: iterative-inequality},  the numerical error $\mathbf{NE}$ is considered as the discrete Lyapunov function, which can be rewritten as 
\begin{equation}
\label{eqn: mono-admm-lyapunov}
\mathcal{E}(k) = \frac{1}{2s} \left\|G(y_{k+1} - y_{k})\right\|^2 + \frac{s}{2} \left\| \lambda_{k+1} - \lambda_k \right\|^2.
\end{equation}
Here, we can calculate its iterative difference as 
\begin{align}
\mathcal{E}(k+1)  - \mathcal{E}(k) & = \underbrace{\frac{1}{s}  \big\langle G(y_{k+2} - y_{k+1}) - G(y_{k+1} - y_{k}), G(y_{k+2} - y_{k+1}) \big\rangle  }_{\mathbf{I}} \nonumber \\
                      & \mathrel{\phantom{=}} + \underbrace{s \big\langle (\lambda_{k+2} - \lambda_{k+1}) - (\lambda_{k+1} - \lambda_{k}), \lambda_{k+2} - \lambda_{k+1} \big\rangle}_{\mathbf{II}}   \label{eqn: mono-admm-lyapunov-iterative} \\
                       & \mathrel{\phantom{=}} - \underbrace{ \frac{1}{2s} \big\| G(y_{k+2} - y_{k+1}) - G(y_{k+1} - y_{k}) \big\|^2  - \frac{s}{2} \big\| (\lambda_{k+2} - \lambda_{k+1}) - (\lambda_{k+1} - \lambda_{k}) \big\|^2 }_{\mathbf{NE}}.  \nonumber 
\end{align}
By utilizing the inequality~\eqref{eqn: mono-admm-3}, we can deduce
\begin{equation}
\label{eqn: mono-admm-lyapunov-iterative-I-II}
\mathbf{I} + \mathbf{II} \leq  \big\langle G(y_{k+2} - y_{k+1}) - G(y_{k+1} - y_{k}), (\lambda_{k+2} - \lambda_{k+1}) - (\lambda_{k+1} - \lambda_{k}) \big\rangle. 
\end{equation}
Combining~\eqref{eqn: mono-admm-lyapunov-iterative} and~\eqref{eqn: mono-admm-lyapunov-iterative-I-II}, we conclude that the iterative sequence satisfies $\mathcal{E}(k+1)  - \mathcal{E}(k) \leq 0$. Therefore, we can enhance~\Cref{thm: admm-convegence-rate} to include the convergence rate with respect to the last iterate. We state this rigorously as follows. 
\begin{theorem}
\label{thm: admm-convegence-rate-mono}
Under the assumption of~\Cref{lem: admm}, the iterative sequence $\{(x_N,y_N; \lambda_N)\}_{N=0}^{\infty}$ generated by~\texttt{ADMM}~\eqref{eqn: admm-x} --- \eqref{eqn: admm-lambda} converges to $(x^{\star},y^{\star};\lambda^{\star})$ with the following rate as
\begin{equation}
\label{eqn:admm-convegence-rate-mono}
\left\|G(y_{N+1} - y_{N})\right\|^2 + s^2 \left\| \lambda_{N+1} - \lambda_{N} \right\|^2  \leq \frac{\left\|G(y_0 - y^{\star})\right\|^2 + s^2 \left\| \lambda_0 - \lambda^{\star} \right\|^2}{N+1}, 
\end{equation}
with any initial $(x_0, y_0; \lambda_0)$.
\end{theorem}

%
%
%
%
%
%
%
%
%
%
%
%
%
%
%
%
%

\section{A general form of ADMM}
\label{sec: general-admm}

In~\Cref{sec: admm} and~\Cref{sec: monotonicity}, we make use of Lyapunov functions, particularly~\eqref{eqn: lypunov-admm-star} and~\eqref{eqn: mono-admm-lyapunov}, to derive the convergence rates of the discrete~\texttt{ADMM} algorithm~\eqref{eqn: admm-x} --- \eqref{eqn: admm-lambda}. However,  it is important to note that the convergence behavior described in these sections does not involve the first iterative sequence $\{x_k\}_{k=0}^{\infty}$ for the general convex functions.  Additionally, as mentioned in~\citep[Section 6.4.1]{boyd2011distributed}, it is still necessary to assume that the matrix $A^{\top}A + sF^{\top}F$ is invertible for the generalized~\texttt{Lasso}. 

In~\citep{he2012convergence, he2015non}, these problems mentioned were effectively resolved by introducing a priori parameter $r$, which is chosen to be greater than the maximum eigenvalue of the matrix $F^{\top}F$, denoted as $\max\lambda(F^{T}F) = \|F^{T}F\|$. By incorporating this parameter, the first iteration of~\texttt{ADMM}~\eqref{eqn: admm-x} is enhanced as a general version, expressed as  
\begin{equation}
\label{eqn: admm-x-extension}
x_{k+1} = \argmin_{x \in \mathbb{R}^d} \left\{ f(x) + \frac{1}{2s} \left( \|Fx + Gy_k - h + s \lambda_k\|^2 + r\|x - x_{k}\|^2 - \|F(x - x_k)\|^2\right) \right\}.
\end{equation}
Combining with the third iteration~\eqref{eqn: admm-lambda}, the general first iteration~\eqref{eqn: admm-x-extension} can be explicitly expressed as
\begin{equation}
\label{eqn: admm-x-extension-explicit}
\partial f(x_{k+1}) - F^{\top}G\cdot \frac{y_{k+1} - y_{k}}{s} + F^{\top} \lambda_{k+1} + (rI_{d_1 \times d_1} - F^{T}F) \cdot \frac{x_{k+1} - x_{k}}{s} \ni 0.
\end{equation}
It is worth noting that the general first iteration~\eqref{eqn: admm-x-extension-explicit} only incorporates a new term $(rI_{d_1 \times d_1} - F^{T}F)(x_{k+1} - x_{k})/s$, which is added to the first iteration~\eqref{eqn: admm-nonsmooth-explicit-x-1}.  

Next, we demonstrate that our theoretical analysis, which is based on the system of high-resolution ODEs and the numerical error resulting from the implicit scheme, is also applicable to the general~\texttt{ADMM},~\eqref{eqn: admm-x-extension},~\eqref{eqn: admm-y} and~\eqref{eqn: admm-lambda}. Moreover, it follows a similar process from the continuous high-resolution ODEs to the general~\texttt{ADMM}, as described in~\Cref{sec: ode} and~\Cref{sec: admm}. For the general first iteration~\eqref{eqn: admm-x-extension-explicit}, we write down its corresponding ODE as
\begin{equation}
\label{eqn: admm-x-extension-ode}
\nabla f(X) - F^{\top}G \dot{Y} + F^{\top}\Lambda + (rI_{d_1 \times d_1} - F^{\top}F) \dot{X} = 0.
\end{equation}
Consistently, the Lyapunov function~\eqref{eqn: lyapunov-ode} is modified and takes the following form as
\begin{equation}
\label{eqn: admm-x-extension-ode-lyapunov}
\mathcal{E}(t) = \frac{r}{2s}\|X - x^{\star}\|^2 - \frac{1}{2s} \|F(X - x^{\star})\|^2 + \frac{1}{2s} \left\|G(Y - y^{\star})\right\|^2 + \frac{s}{2} \left\| \Lambda - \lambda^{\star} \right\|^2.
\end{equation}
With the help of~\Cref{thm: subgrad-unique} and~\Cref{thm: saddle-equivalence}, it is straightforward to show that the time derivative is no more than zero. Furthermore, the discrete Lyapunov function~\eqref{eqn: lypunov-admm-star} is modified as
\begin{equation}
\label{eqn: admm-x-extension-discrete-lyapunov}
\mathcal{E}(k) = \frac{r}{2s}\|x_k - x^{\star}\|^2 - \frac{1}{2s} \|F(x_k - x^{\star})\|^2 + \frac{1}{2s} \left\|G(y_k - y^{\star})\right\|^2 + \frac{s}{2} \left\| \lambda_k - \lambda^{\star} \right\|^2.
\end{equation}
Following a similar process of calculation, we can derive that the iterative difference is bounded by the numerical error as 
\begin{equation}
\label{eqn: admm-x-extension-discrete-ne}
\mathbf{NE}= -  \frac{r}{2s}\|x_{k+1} - x_{k}\|^2 + \frac{1}{2s} \|F(x_{k+1} - x_{k})\|^2  - \frac{1}{2s} \left\|G(y_{k+1} - y_{k})\right\|^2 - \frac{s}{2} \left\| \lambda_{k+1} - \lambda_{k} \right\|^2.
\end{equation}
Therefore, we can enhance~\Cref{thm: admm-convegence-rate} as the following statement.  
\begin{theorem}
\label{thm: admm-convegence-rate-extension}
Under the assumption of~\Cref{lem: admm}, the iterative sequence $\{(x_k,y_k; \lambda_k)\}_{k=0}^{\infty}$ generated by the general~\texttt{ADMM},~\eqref{eqn: admm-x-extension},~\eqref{eqn: admm-y} and \eqref{eqn: admm-lambda}, converges to $(x^{\star},y^{\star};\lambda^{\star})$ with the following rate as claimed in~\eqref{eqn:admm-convegence-rate}. In addition, the iterative sequence $\{(x_k,y_k; \lambda_k)\}_{k=0}^{\infty}$ also obeys
\begin{equation}
\label{eqn:admm-convegence-rate-extension}
\left\|x_{k+1} - x_{k}\right\|^2 \leq \frac{r\left\|x_0 - x^{\star}\right\|^2 + \left\|G(y_0 - y^{\star})\right\|^2 + s^2 \left\| \lambda_0 - \lambda^{\star} \right\|^2}{(N+1)(r - \|F^{\top}F\|)}.
\end{equation}
\end{theorem}

Finally, we demonstrate the monotonicity of the numerical error~\eqref{eqn: admm-x-extension-discrete-ne}, which aligns with the process as described in~\Cref{sec: monotonicity}. We start by considering the continuous system of high-resolution differential equations~\eqref{eqn: admm-x-extension-ode},~\eqref{eqn: admm-ode-y} and~\eqref{eqn: admm-ode-lambda}. By taking the continuous version of the numerical error~\eqref{eqn: admm-x-extension-discrete-ne} as the Lyapunov function, denoted as 
\begin{equation}
\label{eqn: admm-x-extension-ne-Lyapunov}
\mathcal{E}(t) = \frac{rs}{2}\|\dot{X}\|^2 - \frac{s}{2} \|F\dot{X}\|^2 + \frac{s}{2} \|G \dot{Y} \|^2 + \frac{s^3}{2} \| \dot{\Lambda}\|^2
\end{equation}
we can calculate its time derivative. Furthermore, we can determine that it can be probably counterbalanced by the numerical error resulting from the implicit scheme. Then, we implement the same process into the discrete Lyapunov function
\begin{equation}
\label{eqn: admm-x-extension-discrete-ne-Lyapunov}
\mathcal{E}(k) = \frac{r}{2s}\|x_{k+1} - x_{k}\|^2 - \frac{1}{2s} \|F(x_{k+1} - x_{k})\|^2 + \frac{1}{2s} \left\|G(y_{k+1} - y_{k})\right\|^2 + \frac{s}{2} \left\| \lambda_{k+1} - \lambda_{k} \right\|^2.
\end{equation}
By taking some basic analysis, it can be observed that the numerical error causes the iterative difference to decrease monotonically. Therefore, we can enhance~\Cref{thm: admm-convegence-rate-mono} as the following statement. 
\begin{theorem}
\label{thm: admm-convegence-rate-mono-extension}
Under the same conditions with~\Cref{lem: admm}, the iterative sequence $\{(x_N,y_N; \lambda_N)\}_{N=0}^{\infty}$ generated by the general~\texttt{ADMM}, ~\eqref{eqn: admm-x-extension},~\eqref{eqn: admm-y} and \eqref{eqn: admm-lambda}, converges to $(x^{\star},y^{\star};\lambda^{\star})$ with the rate as claimed in~\eqref{eqn:admm-convegence-rate-mono}. In addition, the iterative sequence also obeys
\begin{equation}
\label{eqn:admm-convegence-rate-mono-extension}
\left\|x_{N+1} - x_{N}\right\|^2  \leq \frac{r\left\|x_0 - x^{\star}\right\|^2+ \left\|G(y_0 - y^{\star})\right\|^2 + s^2 \left\| \lambda_0 - \lambda^{\star} \right\|^2}{(N+1)(r - \|F^{\top}F\|)}, 
\end{equation}
with any initial $(x_0, y_0; \lambda_0)$.
\end{theorem}

%
%



\section{Conclusion and discussion}
\label{sec: conclu}

In this paper, we utilize dimensional analysis, as previously employed in~\citep{shi2022understanding, shi2020learning, shi2021hyperparameters}, to derive a system of high-resolution ODEs for the discrete~\texttt{ADMM} algorithm. The system of high-resolution ODEs captures an important feature of~\texttt{ADMM}, called $\lambda$-correction, which causes the trajectory of~\texttt{ADMM} to deviate from the constrained hyperplane. Technically, we utilize Lyapunov analysis to investigate the convergence behavior from the continuous system to the discrete~\texttt{ADMM} algorithm. Through this analysis, we identify that the essential factor that leads to the convergence rate and monotonicity in the discrete case is the numerical error resulting from the implicit scheme. In addition, we further discover that if one component of the objective function is assumed to be strongly convex,  the iterative average of~\texttt{ADMM} converges strongly with a rate $O(1/N)$, where $N$ is the number of iterations.

Understanding and analyzing the discrete~\texttt{ADMM} algorithm through the system of high-resolution ODEs opens up several exciting avenues for further investigation.  It would be intriguing to explore the convergence of~\texttt{ADMM} with various types of norms and rates. Additionally, considering the lower bound $O(1/N^2)$ for the convex optimization problems, as mentioned in~\citep{nemirovski1983problem}, there is significant potential to enhance and accelerate the convergence rates of~\texttt{ADMM}. Furthermore, it is interesting to note that the constrained minimization problem is equivalent to the unconstrained minimax problem, as stated in~\Cref{thm: saddle-equivalence}. In parallel to the~\texttt{ADMM} algorithm for the minimization form,  there is another algorithm for the minimax form, called the~\textit{primal-dual hybrid gradient} algorithm (\texttt{PDHG})~\citep{chambolle2011first, pock2011diagonal, he2012convergence1, chambolle2016introduction}. It would also be fascinating to uncover the mechanism of its iterative convergence using the techniques of ODEs and numerical analysis. Additionally, designing fast algorithms to accelerate~\texttt{PDHG} would be a valuable endeavor. 


%
%
%
%
%
%
%
%


\section*{Acknowledgements}
Bowen Li was partially supported by the Hua Loo-Keng scholarship of CAS.  Bin Shi was partially supported by Grant No.YSBR-034 of CAS.

\bibliographystyle{abbrvnat}
\bibliography{sigproc}

\begin{thebibliography}{36}
\providecommand{\natexlab}[1]{#1}
\providecommand{\url}[1]{\texttt{#1}}
\expandafter\ifx\csname urlstyle\endcsname\relax
  \providecommand{\doi}[1]{doi: #1}\else
  \providecommand{\doi}{doi: \begingroup \urlstyle{rm}\Url}\fi

\bibitem[Adams et~al.(2013)Adams, Dirr, Peletier, and Zimmer]{adams2013large}
S.~Adams, N.~Dirr, M.~Peletier, and J.~Zimmer.
\newblock Large deviations and gradient flows.
\newblock \emph{Philosophical Transactions of the Royal Society A:
  Mathematical, Physical and Engineering Sciences}, 371\penalty0 (20120341),
  2013.

\bibitem[Boyd et~al.(2011)Boyd, Parikh, Chu, Peleato, Eckstein,
  et~al.]{boyd2011distributed}
S.~Boyd, N.~Parikh, E.~Chu, B.~Peleato, J.~Eckstein, et~al.
\newblock Distributed optimization and statistical learning via the alternating
  direction method of multipliers.
\newblock \emph{Foundations and Trends{\textregistered} in Machine learning},
  3\penalty0 (1):\penalty0 1--122, 2011.

\bibitem[Boyd and Vandenberghe(2004)]{boyd2004convex}
S.~P. Boyd and L.~Vandenberghe.
\newblock \emph{Convex optimization}.
\newblock Cambridge University Press, 2004.

\bibitem[Bruckstein et~al.(2009)Bruckstein, Donoho, and
  Elad]{bruckstein2009sparse}
A.~M. Bruckstein, D.~L. Donoho, and M.~Elad.
\newblock From sparse solutions of systems of equations to sparse modeling of
  signals and images.
\newblock \emph{SIAM review}, 51\penalty0 (1):\penalty0 34--81, 2009.

\bibitem[Chambolle and Pock(2011)]{chambolle2011first}
A.~Chambolle and T.~Pock.
\newblock A first-order primal-dual algorithm for convex problems with
  applications to imaging.
\newblock \emph{Journal of mathematical imaging and vision}, 40:\penalty0
  120--145, 2011.

\bibitem[Chambolle and Pock(2016)]{chambolle2016introduction}
A.~Chambolle and T.~Pock.
\newblock An introduction to continuous optimization for imaging.
\newblock \emph{Acta Numerica}, 25:\penalty0 161--319, 2016.

\bibitem[Chen et~al.(2022{\natexlab{a}})Chen, Shi, and Yuan]{chen2022gradient}
S.~Chen, B.~Shi, and Y.-x. Yuan.
\newblock Gradient norm minimization of {N}esterov acceleration: $
  \mathrm{o}(1/k^3) $.
\newblock \emph{arXiv preprint arXiv:2209.08862}, 2022{\natexlab{a}}.

\bibitem[Chen et~al.(2022{\natexlab{b}})Chen, Shi, and
  Yuan]{chen2022revisiting}
S.~Chen, B.~Shi, and Y.-x. Yuan.
\newblock Revisiting the acceleration phenomenon via high-resolution
  differential equations.
\newblock \emph{arXiv preprint arXiv:2212.05700}, 2022{\natexlab{b}}.

\bibitem[Chen et~al.(2023)Chen, Shi, and Yuan]{chen2023underdamped}
S.~Chen, B.~Shi, and Y.-x. Yuan.
\newblock On underdamped {N}esterov's acceleration.
\newblock \emph{arXiv preprint arXiv:2304.14642}, 2023.

\bibitem[Daubechies et~al.(2004)Daubechies, Defrise, and
  De~Mol]{daubechies2004iterative}
I.~Daubechies, M.~Defrise, and C.~De~Mol.
\newblock An iterative thresholding algorithm for linear inverse problems with
  a sparsity constraint.
\newblock \emph{Communications on Pure and Applied Mathematics: A Journal
  Issued by the Courant Institute of Mathematical Sciences}, 57\penalty0
  (11):\penalty0 1413--1457, 2004.

\bibitem[Franca et~al.(2018)Franca, Robinson, and Vidal]{franca2018admm}
G.~Franca, D.~Robinson, and R.~Vidal.
\newblock {ADMM} and accelerated {ADMM} as continuous dynamical systems.
\newblock In \emph{International Conference on Machine Learning}, pages
  1559--1567. PMLR, 2018.

\bibitem[Gabay and Mercier(1976)]{gabay1976dual}
D.~Gabay and B.~Mercier.
\newblock A dual algorithm for the solution of nonlinear variational problems
  via finite element approximation.
\newblock \emph{Computers \& mathematics with applications}, 2\penalty0
  (1):\penalty0 17--40, 1976.

\bibitem[Glowinski and Marroco(1975)]{glowinski1975approximation}
R.~Glowinski and A.~Marroco.
\newblock Sur l'approximation, par {\'e}l{\'e}ments finis d'ordre un, et la
  r{\'e}solution, par p{\'e}nalisation-dualit{\'e} d'une classe de
  probl{\`e}mes de dirichlet non lin{\'e}aires.
\newblock \emph{Revue fran{\c{c}}aise d'automatique, informatique, recherche
  op{\'e}rationnelle. Analyse num{\'e}rique}, 9\penalty0 (R2):\penalty0 41--76,
  1975.

\bibitem[Hastie et~al.(2009)Hastie, Tibshirani, Friedman, and
  Friedman]{hastie2009elements}
T.~Hastie, R.~Tibshirani, J.~H. Friedman, and J.~H. Friedman.
\newblock \emph{The elements of statistical learning: data mining, inference,
  and prediction}, volume~2.
\newblock Springer, 2009.

\bibitem[He and Yuan(2012{\natexlab{a}})]{he2012convergence}
B.~He and X.~Yuan.
\newblock On the $\mathrm{O}(1/n)$ convergence rate of the
  {D}ouglas--{R}achford alternating direction method.
\newblock \emph{SIAM Journal on Numerical Analysis}, 50\penalty0 (2):\penalty0
  700--709, 2012{\natexlab{a}}.

\bibitem[He and Yuan(2012{\natexlab{b}})]{he2012convergence1}
B.~He and X.~Yuan.
\newblock Convergence analysis of primal-dual algorithms for a saddle-point
  problem: from contraction perspective.
\newblock \emph{SIAM Journal on Imaging Sciences}, 5\penalty0 (1):\penalty0
  119--149, 2012{\natexlab{b}}.

\bibitem[He and Yuan(2015)]{he2015non}
B.~He and X.~Yuan.
\newblock On non-ergodic convergence rate of {D}ouglas--{R}achford alternating
  direction method of multipliers.
\newblock \emph{Numerische Mathematik}, 130\penalty0 (3):\penalty0 567--577,
  2015.

\bibitem[Kim et~al.(2009)Kim, Koh, Boyd, and Gorinevsky]{kim2009ell_1}
S.-J. Kim, K.~Koh, S.~Boyd, and D.~Gorinevsky.
\newblock $\ell_1$ trend filtering.
\newblock \emph{SIAM review}, 51\penalty0 (2):\penalty0 339--360, 2009.

\bibitem[Li et~al.(2022{\natexlab{a}})Li, Shi, and Yuan]{li2022linear}
B.~Li, B.~Shi, and Y.-x. Yuan.
\newblock Linear convergence of {ISTA} and {FISTA}.
\newblock \emph{arXiv preprint arXiv:2212.06319}, 2022{\natexlab{a}}.

\bibitem[Li et~al.(2022{\natexlab{b}})Li, Shi, and Yuan]{li2022proximal}
B.~Li, B.~Shi, and Y.-x. Yuan.
\newblock Proximal subgradient norm minimization of {ISTA} and {FISTA}.
\newblock \emph{arXiv preprint arXiv:2211.01610}, 2022{\natexlab{b}}.

\bibitem[Li et~al.(2023)Li, Shi, and Yuan]{li2023linear}
B.~Li, B.~Shi, and Y.-x. Yuan.
\newblock Linear convergence of {N}esterov-1983 with the strong convexity.
\newblock \emph{arXiv preprint arXiv:2306.09694}, 2023.

\bibitem[Nemirovsky and Yudin(1983)]{nemirovski1983problem}
A.~S. Nemirovsky and D.~B. Yudin.
\newblock \emph{Problem Complexity and Method Efficiency in Optimization}.
\newblock Wiley-Interscience,NewYork, 1983.

\bibitem[Nesterov(1998)]{nesterov1998introductory}
Y.~Nesterov.
\newblock \emph{Introductory Lectures on Convex Optimization: A Basic Course},
  volume~87.
\newblock Springer Science \& Business Media, 1998.

\bibitem[Nesterov(1983)]{nesterov1983method}
Y.~E. Nesterov.
\newblock A method of solving a convex programming problem with convergence
  rate $o\bigl( 1/k^2 \bigr)$.
\newblock \emph{Doklady Akademii Nauk}, 269\penalty0 (3):\penalty0 543--547,
  1983.

\bibitem[Pock and Chambolle(2011)]{pock2011diagonal}
T.~Pock and A.~Chambolle.
\newblock Diagonal preconditioning for first order primal-dual algorithms in
  convex optimization.
\newblock In \emph{2011 International Conference on Computer Vision}, pages
  1762--1769. IEEE, 2011.

\bibitem[Rockafellar(1970)]{rockafellar1970convex}
R.~T. Rockafellar.
\newblock \emph{Convex Analysis}, volume~18.
\newblock Princeton University Press, 1970.

\bibitem[Rockafellar and Wets(2009)]{rockafellar2009variational}
R.~T. Rockafellar and R.~J.-B. Wets.
\newblock \emph{Variational analysis}, volume 317.
\newblock Springer Science \& Business Media, 2009.

\bibitem[Rudin et~al.(1992)Rudin, Osher, and Fatemi]{rudin1992nonlinear}
L.~I. Rudin, S.~Osher, and E.~Fatemi.
\newblock Nonlinear total variation based noise removal algorithms.
\newblock \emph{Physica D: nonlinear phenomena}, 60\penalty0 (1-4):\penalty0
  259--268, 1992.

\bibitem[Shi(2021)]{shi2021hyperparameters}
B.~Shi.
\newblock On the hyperparameters in stochastic gradient descent with momentum.
\newblock \emph{arXiv preprint arXiv:2108.03947}, 2021.

\bibitem[Shi et~al.(2019)Shi, Du, Su, and Jordan]{shi2019acceleration}
B.~Shi, S.~S. Du, W.~Su, and M.~I. Jordan.
\newblock Acceleration via symplectic discretization of high-resolution
  differential equations.
\newblock \emph{Advances in Neural Information Processing Systems}, 32, 2019.

\bibitem[Shi et~al.(2020)Shi, Su, and Jordan]{shi2020learning}
B.~Shi, W.~J. Su, and M.~I. Jordan.
\newblock On learning rates and schr$\backslash$" odinger operators.
\newblock \emph{arXiv preprint arXiv:2004.06977}, 2020.

\bibitem[Shi et~al.(2022)Shi, Du, Jordan, and Su]{shi2022understanding}
B.~Shi, S.~S. Du, M.~I. Jordan, and W.~J. Su.
\newblock Understanding the acceleration phenomenon via high-resolution
  differential equations.
\newblock \emph{Mathematical Programming}, 195\penalty0 (1-2):\penalty0
  79--148, 2022.

\bibitem[Su et~al.(2016)Su, Boyd, and Candes]{su2016differential}
W.~Su, S.~Boyd, and E.~J. Candes.
\newblock A differential equation for modeling {N}esterov's accelerated
  gradient method: Theory and insights.
\newblock \emph{Journal of Machine Learning Research}, 17\penalty0
  (153):\penalty0 5312–5354, 2016.

\bibitem[Tibshirani(1996)]{tibshirani1996regression}
R.~Tibshirani.
\newblock Regression shrinkage and selection via the lasso.
\newblock \emph{Journal of the Royal Statistical Society Series B: Statistical
  Methodology}, 58\penalty0 (1):\penalty0 267--288, 1996.

\bibitem[Wilson et~al.(2021)Wilson, Recht, and Jordan]{wilson2021lyapunov}
A.~C. Wilson, B.~Recht, and M.~I. Jordan.
\newblock A {L}yapunov analysis of accelerated methods in optimization.
\newblock \emph{Journal of Machine Learning Research}, 22\penalty0
  (113):\penalty0 5040--5073, 2021.

\bibitem[Wooldridge(2012)]{wooldridge2012introductory}
J.~M. Wooldridge.
\newblock \emph{Introductory Econometrics: A Modern Approach}.
\newblock South-Western, Cengage Learning, fifth edition, 2012.

\end{thebibliography}
\end{document}